\documentclass[preprint,12pt]{elsarticle}

\usepackage{amssymb}
\usepackage{amsmath}

\usepackage[table,xcdraw]{xcolor}
\usepackage{bm}
\usepackage{algorithm2e}
\usepackage{subcaption}
\usepackage{multirow}
\usepackage{diagbox}
\usepackage{colortbl}
\usepackage{xcolor}  
\usepackage{makecell}
\usepackage{float}
\usepackage{subcaption}
\usepackage{hyperref}
\usepackage{cleveref}
\usepackage{enumerate}
\graphicspath{{media/}} 
\usepackage[inline]{enumitem}
\usepackage{todonotes}

\journal{JCP}
\newcommand{\xnet}{IV-Net}
\newcommand{\nm}{\mathbin{\phantom{-}}}
\newcommand{\nmm}{\mathbin{\phantom{0}}}

\newif\iffinal

\begin{document}

\begin{frontmatter}



\title{IV-Net: A neural  network for elliptic PDEs with random and highly varying coefficients}

%

\author[label1]{Shan Zhong}

\author[label1,label2]{George Biros}
\affiliation[label1]{organization={Oden Institute for Computational Science and Engineering, The University of Texas at Austin},
            addressline={201 E 24th St},
            city={Austin},
            postcode={78712},
            state={TX},
            country={USA}}
\affiliation[label2]{organization={Walker Department of Mechanical Engineering, The University of Texas at Austin},
            addressline={301 E Dean Keeton St},
            city={Austin},
            postcode={78712},
            state={TX},
            country={USA}}


\begin{abstract}
We introduce a novel neural operator architecture designed to approximate solutions of linear elliptic partial differential equations with high-contrast, spatially varying coefficients. The network, termed the \textbf{I}terated \textbf{V}-shaped \textbf{Net} (\textbf{\emph{\xnet}}), realizes a mapping from the input coefficients and right-hand side to the corresponding solution field. The architecture of \xnet\ is informed by, and closely resembles, a V-cycle multigrid solver. 

The \xnet\ model is parameterized via convolutional layers defined in the physical domain. For coercive problems with highly heterogeneous coefficients, the proposed network exhibits superior performance relative to a proper orthogonal decomposition (POD) approach and several existing neural operator architectures.  For low-frequency oscillatory Helmholtz problems with smooth coefficients, its performance is similar to  that of a Fourier neural operator. 

We analyze the approximation error and convergence behavior of \xnet, its data efficiency, and its dependence on the underlying discretization mesh. Furthermore, we demonstrate the practical effectiveness of the architecture through a series of numerical experiments, including applications to uncertainty quantification, inverse problems, and prediction of quantities of interest.
\end{abstract}

\end{frontmatter}


\section{Introduction}
We are interested in constructing a neural network surrogate for the solution of the following Partial Differential Equation (PDE) with variable coefficients $\eta, \eta$ and $a$:
\begin{equation}
\label{e:proto}
(2\pi\kappa)^2(1+\eta(x \xi))u(x) + \beta \nabla(a(x)\nabla u(x)) = f(x, \lambda), \quad x\in 
\Omega :=(0,1)^2, 
\end{equation}
where $\kappa$ is a non-negative scaling constant, $\xi$ and $\lambda$ denote the random variables that generate the fields $\eta \in L^{\infty}((0,1)^2;\mathbb{R}_+)$ and $f \in L^2((0,1)^2;\mathbb{R})$, $\beta \in \left\{-1,1\right\}$, $a$ is also a random field to be specified later.
Our goal is to approximate the mapping $(\eta, a, f)\rightarrow u$.\par
To differentiate between different problem instances and with a slight abuse of conventional definitions, we refer to the case $\kappa >0, \beta = -1, a = 1$ as the "Poisson problem", $\kappa =0, \beta = -1, a >0$ as the "Darcy problem", and $\kappa >0, \beta = 1, a = 1$ as the "Helmholtz problem".%
\footnote{We include the Helmholtz problem as a proxy for the complex Helmholtz operator. Also note that for simplicity, we use Neumann conditions as opposed to absorbing ones.}
The Poisson and Darcy problems are coercive. We use Neumann boundary conditions for the Poisson and Helmholtz problems, and Dirichlet conditions for the Darcy problem.\par
We only consider problems on a square domain using a regular grid discretization. To generate training and testing data, we discretize \Cref{e:proto} using a standard $2^{\text{nd}}$ order finite difference scheme. Using  bold face letters for the discretized operator and related fields,  \cref{e:proto} reads
\begin{equation}
\label{e:discretized_proto}
\mathbf{A}[\bm{\eta}, \mathbf{a}] \mathbf{u} = \mathbf{f},
\end{equation}
where $ \bm{\eta}, \mathbf{f}, \mathbf{a} \in \mathbb{R}^n$ are discretized fields, $\mathbf{u}\in \mathbb{R}^n$ is the solution, $\mathbf{A} \in \mathbb{R}^{n\times n}$ is a sparse matrix depending on $\bm{\eta}$ and $\mathbf{A}$, and $n$ is the number of unknowns. \xnet\ aims to approximate $\mathbf{u}=\mathbf{A}[\bm{\eta}, \mathbf{a}]^{-1} \mathbf{f}$. We also include a brief discussion on using \xnet\ for the inverse problem $\{\mathbf{f}, \mathbf{u}\} \rightarrow \bm{\eta}.$ 

The surrogate model can be employed in a variety of many-query settings such as uncertainty quantification (UQ), Bayesian inference, and optimal control. These tasks typically require repeated solutions of the underlying PDE for varying parameter values at each stage or inner iteration, which can incur substantial computational cost. Within this framework, the primary motivation for constructing surrogates for \cref{e:proto} is to {\em accelerate} the repeated solution of \cref{e:discretized_proto} for different realizations of $\bm{\eta}$ and $\mathbf{a}$. A further motivation for investigating this problem is to gain deeper insight into the behavior of approximations to neural operator architectures, which are expected to be applicable in more challenging problem settings. 

The reader may rightly question what is meant here by “acceleration.” For coercive problems, domain decomposition and multigrid methods are already known to attain optimal $O(n)$ complexity \cite{chenMetaMgNetMetaMultigrid2020}. In our setting, in which no prior knowledge of a lower-dimensional representation of the input fields is assumed, the surrogate model must also have at least $O(n)$ complexity, as it is necessary to read all coefficients. Consequently, even in the most favorable scenario, our surrogate can at best yield improvements in the multiplicative constants of the overall complexity estimate, rather than in the asymptotic scaling itself. Its practical utility will therefore depend on problem-specific characteristics, implementation details, and the desired target accuracy. In this work, we restrict our attention to demonstrating the feasibility of constructing a surrogate based on a neural network, \emph{without} claiming that such a surrogate surpasses established fast solvers in performance.

There has been significant interest in recent years to develop data driven machine learning methods for the fast, accurate solution of PDEs. 
Earlier approaches, such as physics-informed or physics-constrained learning methods \cite{eDeepRitzMethod2018, raissiPhysicsinformedNeuralNetworks2019, sitzmannImplicitNeuralRepresentations2020},  are designed as alternative PDE discretization methods.  Neural operator methods have been developed to learn the map between the parameter and solution spaces of PDEs, so that a family of PDEs with varying parametric terms can be efficiently solved by one model. Some important works in this area include networks using convolutional encoder-decoder (downsampling-upsampling) layers \cite{zhuBayesianDeepConvolutional2018,khooSolvingParametricPDE2021, raonicCONVOLUTIONALNEURALOPERATORS2023, heMgNOEfficientParameterization2023, winovichConvPDEUQConvolutionalNeural2019}, networks based on universal operator approximation theorem \cite{chenUniversalApproximationNonlinear1995,luLearningNonlinearOperators2021,kontolatiLearningNonlinearOperators2024}, and networks parameterizing the map in spectral spaces, such as Fourier bases \cite{liFourierNeuralOperator2021, rahmanUshapedNeuralOperators2023, wenUFNOEnhancedFourier2022, kovachkiNeuralOperatorLearning}, proper decomposition bases \cite{bhattacharyaModelReductionNeural2021, luComprehensiveFairComparison2022, oleary-roseberryDerivativeInformedNeuralOperator2024}, and multiwavelet bases \cite{guptaMultiwaveletbasedOperatorLearning2021}. The attention mechanism used in \cite{vaswaniAttentionAllYou2017a} has also been explored towards learning operators \cite{caoChooseTransformerFourier2021, liTransformerPartialDifferential2023, yePDEformerFoundationModel2024, kissasLearningOperatorsCoupled2022}.
We discuss selected network architectures in detail in \Cref{s:related}.
\par

{\bf Contributions:}
This work presents an empirical investigation based exclusively on numerical experiments.
In this context, our contributions are twofold. First, we introduce \xnet\, which, as we demonstrate, exhibits characteristics analogous to those of an iterative numerical method.
Second, we propose a set of challenging benchmark test cases and systematically assess the performance of widely used neural network architectures on these new scenarios.
For each test case, we conduct convergence studies with respect to the size of the training dataset, the model capacity, and the resolution of the underlying computational mesh.

To facilitate a comparison with standard model reduction techniques, we contrast the neural operator approximations with a baseline proper orthogonal decomposition (POD) method (\Cref{s:proj}). In addition to POD, we evaluate the following neural operator architectures: U-Net \cite{ronnebergerUNetConvolutionalNetworks2015a}, the Fourier Neural Operator (FNO) \cite{liFourierNeuralOperator2021}, and DeepONet \cite{luLearningNonlinearOperators2021}. For these architectures, we consider multiple network configurations and tune hyperparameters using the same optimization protocol as for \xnet.

Our numerical experiments indicate that, for coercive problems with strongly heterogeneous coefficients, \xnet\ can substantially outperform both POD and the other neural operator architectures. Furthermore, we observe that, although the loss function and mean approximation error can be small, the worst-case error during both training and inference may still be considerable. In \Cref{s:sec_diss}, we discuss several possible extensions and modifications of \xnet. We conclude with an analysis of the limitations of \xnet\ and a broader perspective on its applicability in \Cref{s:conclusion}.

The remainder of the paper is organized as follows. \Cref{s:sec_method} details the equations to be studied. Next, we introduce a proper orthogonal decomposition method, and explain the proposed network. We also discuss related works, and highlight the novelty of our network. In \Cref{s:sec_results}, we present experiments to assess \xnet's accuracy, error convergence rate, data efficiency, mesh adaptivity, performance in uncertainty quantification predictions, and make comparisons with other methods. 
Then, we discuss direct predictions of Quantities of Interest (QoI), application of \xnet\ in inverse problems, and ablation studies of key architectural choices in \Cref{s:sec_diss}. We point out the limitations of our approach in \Cref{s:conclusion}. Additional technical details and experiments are included in the appendix. We summarize the notation in \Cref{t:notation}.\par

\begin{table}[]
\centering
\caption{Table of Notation}
\label{t:notation}
\resizebox{\textwidth}{!}{
\begin{tabular}{ll}
Notation                             & Meaning                                                                     \\ \hline
\multicolumn{1}{|l|}{\bf Problem formulation} & \multicolumn{1}{l|}{}                                                       \\
\multicolumn{1}{|l|}{$\Omega$}              & \multicolumn{1}{l|}{The spatial domain for the PDEs}                        \\
\multicolumn{1}{|l|}{$x$}              & \multicolumn{1}{l|}{Points in the spatial domain}                             \\
\multicolumn{1}{|l|}{$f$}              & \multicolumn{1}{l|}{The source term}                             \\
\multicolumn{1}{|l|}{$\eta, a$}              & \multicolumn{1}{l|}{The varying coefficients in \Cref{e:proto}}                             \\
\multicolumn{1}{|l|}{$\mathbf{f}$}              & \multicolumn{1}{l|}{Discretized source term}                             \\
\multicolumn{1}{|l|}{$\bm{\eta}, \mathbf{a}$}              & \multicolumn{1}{l|}{Discretized varying coefficients}                             \\
\multicolumn{1}{|l|}{$\mathbf{A}$}              & \multicolumn{1}{l|}{Discretized PDE operator}                             \\
\multicolumn{1}{|l|}{$\xi, \lambda$}              & \multicolumn{1}{l|}{The intrinsic random variables for $\eta, f$}                             \\
\multicolumn{1}{|l|}{$c$}              & \multicolumn{1}{l|}{The contrast factor for the Poisson problems}                             \\
\multicolumn{1}{|l|}{$\kappa$}              & \multicolumn{1}{l|}{The wave number in the Helmholtz problems}       \\
\multicolumn{1}{|l|}{$\zeta$}              & \multicolumn{1}{l|}{A probability measure to be sampled from}               \\
\hline

\multicolumn{1}{|l|}{\bf Neural network} & \multicolumn{1}{l|}{}                                                       \\
\multicolumn{1}{|l|}{$V$}              & \multicolumn{1}{l|}{A V-shaped block performing one iterative update}       \\
\multicolumn{1}{|l|}{$\mathcal{P}$}              & \multicolumn{1}{l|}{The prolongation operation}                             \\
\multicolumn{1}{|l|}{$\mathcal{R}$}              & \multicolumn{1}{l|}{The restriction operation}                               \\ 
\multicolumn{1}{|l|}{$C, T$}              & \multicolumn{1}{l|}{Convolution and transposed convolution modules}                             \\
\multicolumn{1}{|l|}{$\sigma$}              & \multicolumn{1}{l|}{The activation function}                             \\
\multicolumn{1}{|l|}{$\mathcal{N}$}              & \multicolumn{1}{l|}{The batch normalization operation}                             \\
\multicolumn{1}{|l|}{$n_d$}              & \multicolumn{1}{l|}{The number of different  dilation parameters}                             \\
\multicolumn{1}{|l|}{$n_b$}              & \multicolumn{1}{l|}{The number of blocks}                             \\
\multicolumn{1}{|l|}{$L$}              & \multicolumn{1}{l|}{The number of different resolution levels}                             \\
\multicolumn{1}{|l|}{$w$}              & \multicolumn{1}{l|}{The multiplicative factor of the number of channels for convolutions}                             \\
\multicolumn{1}{|l|}{$\mathcal{L}$}              & \multicolumn{1}{l|}{The loss function}                             \\\hline
\multicolumn{1}{|l|}{\bf Parameters} & \multicolumn{1}{l|}{}                                                       \\
\multicolumn{1}{|l|}{$n$}              & \multicolumn{1}{l|}{Total number of spatial unknowns for $\bm{\eta}, \mathbf{f}, \mathbf{u}$ (assumed equal)}                             \\
\multicolumn{1}{|l|}{$s$}              & \multicolumn{1}{l|}{The number of training samples}                             \\
\multicolumn{1}{|l|}{$M$}              & \multicolumn{1}{l|}{Total number of weights in a neural network}                             \\
\multicolumn{1}{|l|}{$d_q$}              & \multicolumn{1}{l|}{Dimension of a quantity of interest $q$}                             \\
\hline
\end{tabular}
}
\end{table}


\section{Forward PDE Setup and Discretization}
\label{s:sec_method}
In the Poisson and Helmholtz problems, we assume that $\eta, f, a$ are not arbitrary but  enjoy some, in general unknown,  lower dimensional structure instead.
In the Darcy problem $a$ is infinite dimensional.
In all cases, for the purposes of testing, we constract synthetic fields by considering a stochastic model for $\eta(x), f(x), a(x)$.
In the following $\xi \in \mathbb{R}^{d_{\xi}}, \lambda  \in \mathbb{R}^{d_{\lambda}}$ are random variables that parametrize $\eta = \eta(x;\xi), f = f(x;\lambda)$.
\subsection{Problem Formulation}
\label{s:setup}

\begin{figure}
    \centering
    \begin{subfigure}{1.0\linewidth}
        \centering
        \includegraphics[width=0.9\linewidth]{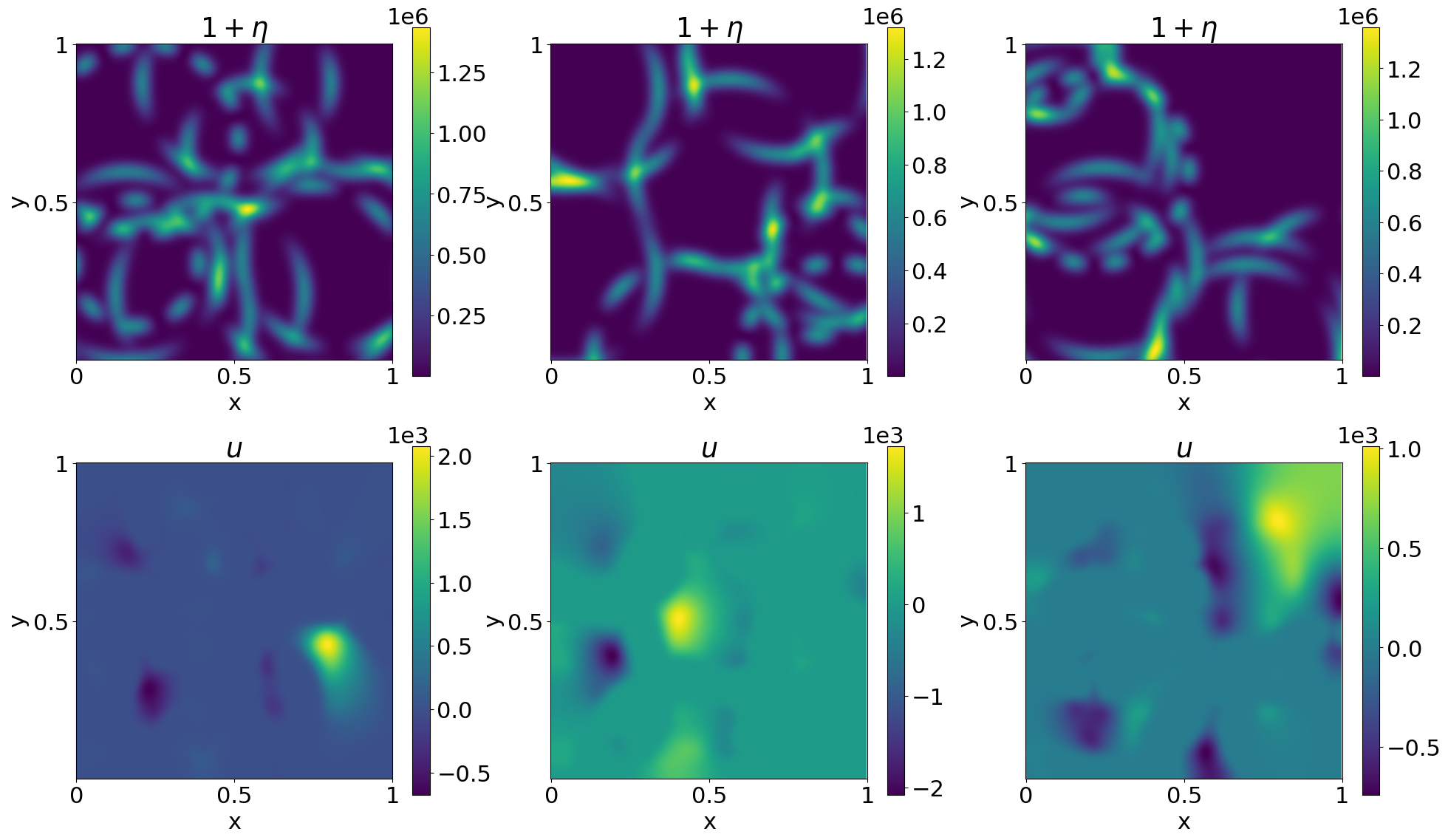}
        \caption{Example coefficients (upper row) and solutions (lower row) of the Poisson problem.}
    \end{subfigure}
    \smallskip
    \begin{subfigure}{1.0\linewidth}
        \centering
        \includegraphics[width=0.9\linewidth]{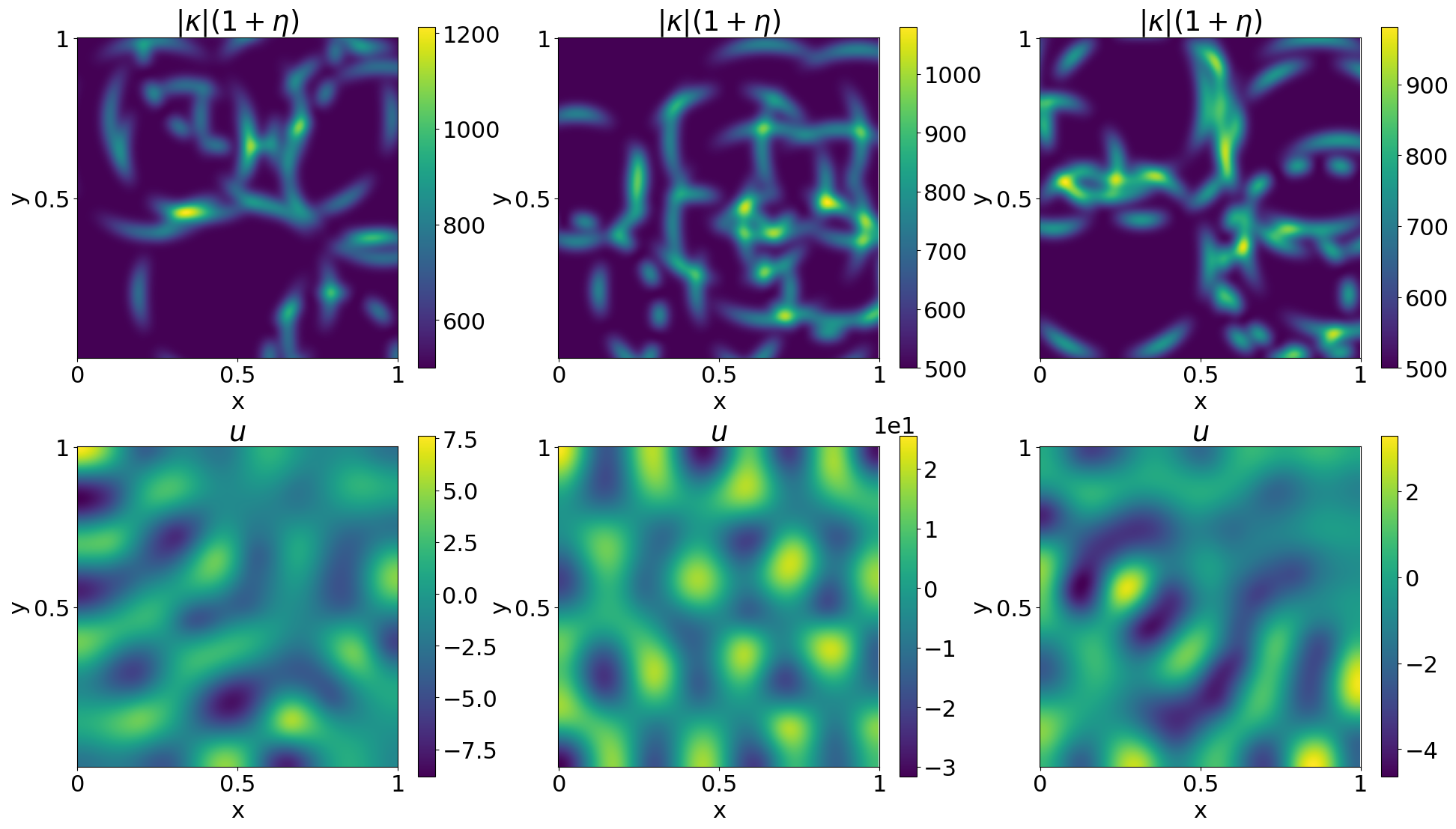}
        \caption{Example coefficients (upper row) and solutions (lower row) of the Helmholtz problem.}
    \end{subfigure}
\end{figure}

\subsubsection{The Poisson Problem}
With $(2\pi\kappa)^2=1$, $\beta=-1$, $a(x)=1$ in \cref{e:proto}, we consider the following PDE with Neumann boundary conditions:
\begin{equation}
\label{e:gh}
\begin{split}
    \left( 1+\eta(x, \xi)\right)u(x)-\Delta u(x) &= f(x, \lambda), \quad x\in(0,1)^2, \\
    \frac{\partial u}{\partial n} &= g. 
\end{split}
\end{equation}
We design $\eta(x)$ to have large contrast and consist of many islands of high contrast, which are known to be challenging for multigrid solvers.
Specifically,  we represent $\eta$ as a sum of ten randomly placed rings with angular modulations to simulate $\eta$ (see \Cref{f:examples}(a) for examples):
\begin{equation}
\label{e:bh}
    \eta(x) = c \sum_{i=1}^{10}{\rm cos^2}\left(2k_i\theta_i(x)\right)
    {\rm exp}\left( -\frac{(\|x-x_i\|-r_i)^2}{2\sigma^2}\right), \quad c=6\times10^5,
\end{equation}
where $x_i$, $r_i$ denote the centers and radii of the rings, and the cosine function is used to generate angular modulations with frequencies controlled by $k_i$. Each $\theta_i(x)$ is the local polar angle defined relative to the $i_{\mathrm{th}}$ ring center. For the following experiments, we assume a fixed right-hand-side $f$, but we also discuss generalization to varying $f$ in \Cref{s:NN}. The randomness of $\eta$ comes from center positions $\xi = \{x_i\}_{i=1}^{d_{\xi}} \in \mathbb{R}^{2d_{\xi}}$, where $d_{\xi}=10$. Each $x_i$ is uniformly distributed over the whole domain. Thus, the intrinsic dimension of $\eta$ is 20, determined by the center positions of the rings. We define the contrast of coefficient$ 1+ \bm{\eta}$ as $\frac{\text{max}(1 +\bm{\eta})}{\text{min}(1+\bm{\eta})}$, which is directly controlled by the scaling factor $c$ in \cref{e:bh}, which is quite large. Other parameters $\{k_i\}$, $\{r_i\}$ are fixed. The examples of $\bm{\eta}$ and the corresponding solutions are shown in \Cref{f:examples}. Such $\bm{\eta}$ resembles coefficients found in  problems in reservoir modeling and flows in fractured porous media.

\subsubsection{The Helmholtz Problem}
We adopt the parametric setting of $\eta$ in the Poisson case, but explore the Helmholtz equation ($\kappa$ positive, $\beta=1$, $a(x)=1$ in \cref{e:proto}):
\begin{equation}
\label{e:darcy}
\begin{split}
    -(2\pi\kappa)^2\left( 1+\eta(x, \xi)\right)u(x)-\Delta u(x) &= f(x), \quad x\in(0,1)^2, \\
    \frac{\partial u}{\partial n} &= g. 
\end{split}
\end{equation}
We use a different scaling constant $c=0.5$ in $\eta$ defined in \cref{e:bh} so that $\|\eta\|_{\infty} \sim 1$. In addition to $\eta$, the difficulty of the problem is also controlled by the wave number $\kappa$. A larger $\kappa$ allows propagation of higher wave numbers and therefore more spatially complex wave oscillation patterns in the solution \cite{yuMultiScaleContextAggregation2016}. For a fixed grid of $256\times256$, we select $(2\pi\kappa)^2=100$ and $500$. This obeys the rule of thumb that at least 10 grid points per wavelength are used for enough resolution 
\cite{babuskaPollutionEffectFEM1997, azulayMultigridAugmentedDeepLearning2023}. 
Examples of the $\eta$ coefficients ($(2\pi\kappa)^2=500$) and the corresponding solutions are shown in \Cref{f:examples}(b).\par

\subsubsection{The Darcy Problem}
\begin{figure} [h!]
    \begin{subfigure}{1.0\linewidth}
        \centering
        \includegraphics[width=0.9\linewidth]{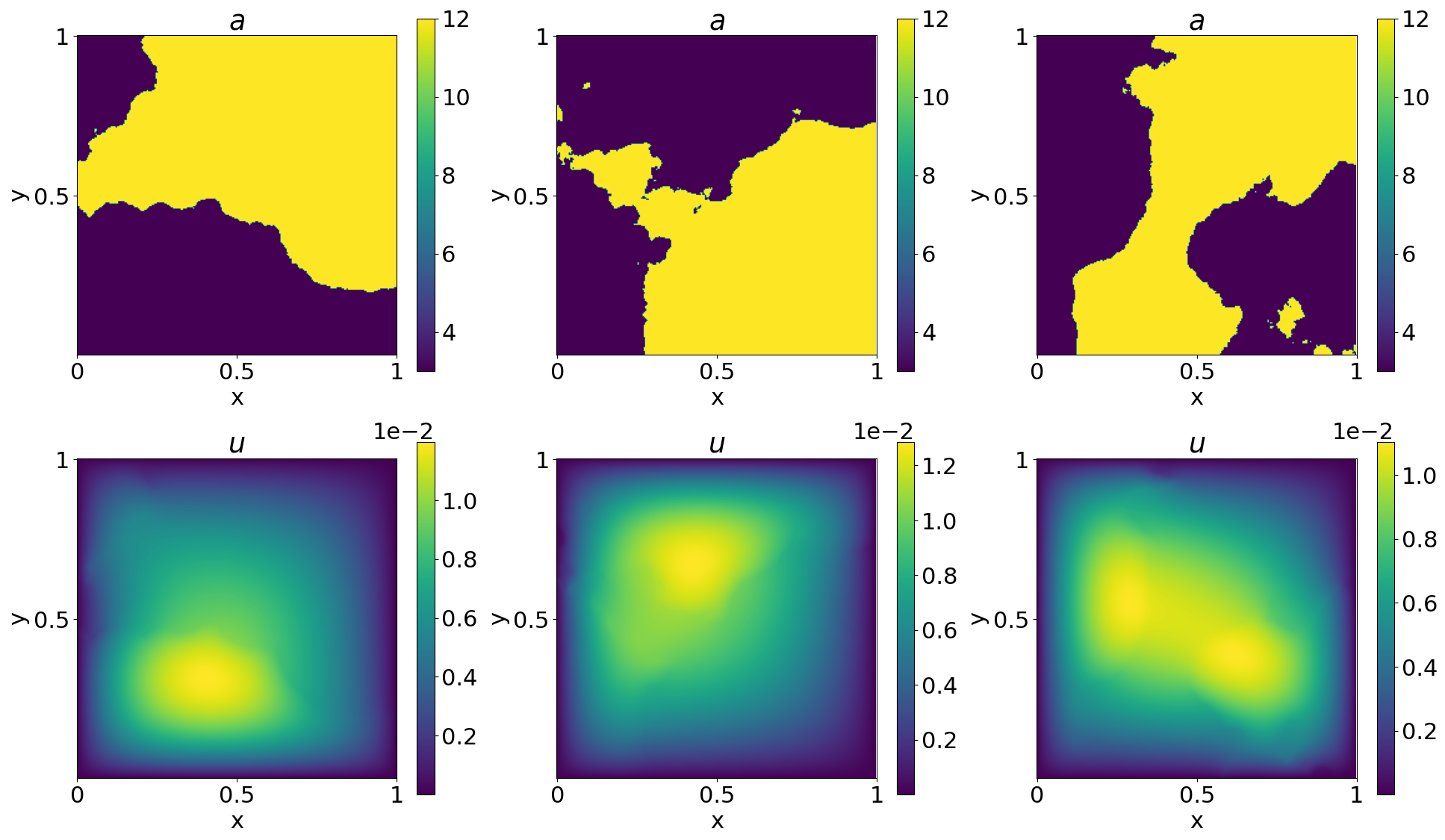}
        \caption{Example coefficients (upper row) and solutions (lower row) of the Darcy problem.}
    \end{subfigure}
    \caption{Example coefficients and solutions for Poisson, Helmholtz, and Darcy problem. A sum of ten randomly placed ring-like patterns simulates material properties in both Poisson and Helmholtz problems. The solutions of the Poisson problem are non-oscillatory, while oscillatory features are observed in Helmholtz cases. For the Darcy problem, the coefficients are binary-valued and consist of several regions of constant values. }
    \label{f:examples}
\end{figure}

The two-dimensional Darcy flow equation is a widely-studied test case in the literature, so we also select this problem here for comparison with Fourier Neural Operator and DeepONet. It describes the single phase, steady-state flow through a random permeability field $a$. It is derived by letting $\kappa=0$ in \cref{e:proto}:
\begin{equation}
\begin{split}
    -\nabla(a(x)\nabla u(x)) &= f(x), \quad x\in (0,1)^2, \\
    u(x) &= 0, \quad x\in \partial (0,1)^2.
\end{split}
\end{equation}
The dataset is generated as described in \cite{liFourierNeuralOperator2021, rahmanUshapedNeuralOperators2023}, where $a \sim \zeta$, $\zeta$ is the push-forward measure of $\mathcal{N}(0, (-\Delta +9I)^{-2})$ under $\psi$, and $f(x)=1$. $\psi(x)$ is defined as:
\begin{equation}
\label{e:phi}
\psi(x)= \begin{cases}3 & \text { if } x<0 \\ 12 & \text { if } x \geq 0\end{cases}.
\end{equation}  
This creates a piecewise-constant binary-valued field $a$ composed of several constant areas (See \Cref{f:examples}(c)). 

\subsection{Proper Orthogonal Decomposition Method}
\label{s:proj}
Let us first discuss the proper orthogonal decomposition method. It will serve as a baseline to be compared with deep learning models, and help us assess the difficulty of each problem. Notice that more sophisticated, manifold-based methods exist, but only a simple method is used here to set up a baseline. For simplicity, we describe the method for the Poisson Problem with $a(x)=1$.
\par
Suppose we have $s$ solutions $\{u_i\}_{i=1}^s$, $u_i \in \mathbb{R}^n$,  computed by sampling $\bm{\eta}$ and solving corresponding equations ($s\ll n$). Let $\mathbf{X}\in \mathbb{R}^{n\times s}$, $\mathbf{X} = (\mathbf{u}_1, \mathbf{u}_2, ..., \mathbf{u}_s)$ be the matrix having these $s$ solutions as columns, and let its singular value decomposition (based on $L_2$ norm) be $\mathbf{X=USV}^{\rm T}$. 
We can truncate the first $r$ columns of $\mathbf{U}$ as $\mathbf{U}_r\in \mathbb{R}^{n\times r}$, and this is the optimal $r$-dimensional subspace for projection \cite{bennerSurveyProjectionBasedModel2015, bhattacharyaModelReductionNeural2021}. The matrix $\mathbf{U_r}\in \mathbb{R}^{n\times r}$ will be used to define a Galerkin scheme where $\mathbf{u} = \mathbf{U}\mathbf{\hat{u_r}}$. 
The discretized linear system is projected onto that subspace,
\begin{equation}\label{e:7}
    \mathbf{U}_r^{\rm T}\mathbf{A}[\bm{\eta}]\mathbf{U}_r \hat{\mathbf{u}}_r = \mathbf{U}_r^{\rm T}\mathbf{f}.
\end{equation}
Denote 
\begin{equation}
    \mathbf{A}_r = \mathbf{U}_r^{\rm T}\mathbf{A}[\bm{\eta}]\mathbf{U}_r, 
    \quad \mathbf{f}_r = \mathbf{U}_r^{\rm T}\mathbf{f}.
\end{equation}
Since we assume $r \ll n$, we can solve for $\hat{\mathbf{u}}_r$ using an LU factorization:
\begin{equation}
    \mathbf{A}_r[\bm{\eta}]\hat{\mathbf{u}}_r = \mathbf{f}_r, \quad \hat{\mathbf{u}}_r = \mathbf{A}_r[\bm{\eta}]^{-1}\mathbf{f}_r.
\end{equation}
\par
Applying $\mathbf{U}_r^{\rm T}$ costs $O(nr)$ operations. Inverting $\mathbf{A}_r$ requires $O(r^3)$ operations. So, assuming $r\ll n$, the overall cost is  $O(nr)$ for one instance of $\bm{\eta}$. 
If $\bm{\eta}$ is a linear combination of random fields such as \Cref{e:bh}, then $\mathbf{A}_r[\bm{\eta}]$ can be precomputed in $O(r^2)$. Otherwise, in general,  $\mathbf{A}_r[\bm{\eta}]$ cannot be precomputed, thus the $O(nr)$ complexity. Let us remark that our POD approach is elementary. More advanced POD methods exist that can better handle parameterized PDEs~\cite{quarteroni2015reduced}. Its purpose is to serve as a reference due to its simple implementation.

\subsection{\xnet}
\label{s:NN}
Let us now turn our attention to the proposed \xnet.

Let's assume that $\mathbf{A}$ is SPD and consider solving the discretized equation $\mathbf{A}[\bm{\eta}]\mathbf{u}=\mathbf{f}$ via steepest descent for 
\begin{equation}
    \min_{\mathbf{u}} \quad \frac{1}{2}\mathbf{u}^{\rm T}\mathbf{A}[\bm{\eta}]\mathbf{u}-\mathbf{u}^{\rm T}\mathbf{f}.
\end{equation}
The steepest descentd update reads
\begin{equation}
\label{e:update}
    \mathbf{u}_{k+1} = \mathbf{u}_{k} - \alpha_k (\mathbf{A}[\bm{\eta}]\mathbf{u}_k-\mathbf{f}),
\end{equation}
where $\alpha_k$ is the step size.
(We remark that this scheme is also related to connection between ODE solvers  and deep networks first proposed in \cite{eProposalMachineLearning2017}.) 
We can construct a network block to compute $\mathbf{u}_{k+1}$: the $\mathbf{u}_k$ term in \cref{e:update} is delivered by the residual connection \cite{heDeepResidualLearning2016, eProposalMachineLearning2017, longPDENetLearningPDEs2018}, and the rest part is approximated by a V-shaped multilevel convolutional structure.
\par
To describe \xnet, we define the following functions: 
\begin{itemize}
    \item $\mathcal{\sigma}: \mathbb{R} \rightarrow \mathbb{R}$, an element-wise nonlinear activation function,
    \item $\mathcal{N}: \mathbb{R}^{c\times n} \rightarrow \mathbb{R}^{c\times n}$, batch normalization \cite{ioffeBatchNormalizationAccelerating2015}; given each channel $x_i\in \mathbb{R}^n$, $i=1,\dots, c$,  it computes $\mathcal{N}(x_i) \in \mathbb{R}^n$,
    $$\mathcal{N}(x_i) = \frac{x_i-\mathrm{E}[x_i]}{\sqrt{\operatorname{Var}[x_i]+\epsilon}} * \gamma_i+\beta_i, \quad i=1,\dots, c,$$
    where $\gamma_i, \beta_i \in \mathbb{R}$ are learnable parameters, $\mathrm{E}[x_i], \operatorname{Var}[x_i] \in \mathbb{R}$ are running estimations based on historical training batches,
    
    \item $\mathcal{C}_{i,j}: \mathcal{R}^{c_{\text{in}}\times n_{\text{in}}} \rightarrow \mathbb{R}^{c_{\text{out}}\times n_{\text{out}}}$, the convolution in the $i_{\mathrm{th}}$ restriction step $\mathcal{R}_i$ with the $j_{\mathrm{th}}$ dilation parameter, and $n_{\text{in}} > n_{\text{out}}$,
    \item $\mathcal{R}_i^j \bm{h} = \sigma(\mathcal{C}_{i,j} \bm{h})$, the restriction operation with the $j_{\mathrm{th}}$ dilation parameter in the $i_{\mathrm{th}}$ restriction step,
    \item $\mathcal{R}_i \bm{h} = \mathcal{N}([\mathcal{R}_i^1 \bm{h}, \mathcal{R}_i^2 \bm{h}, \dots, \mathcal{R}_i^{n_d} \bm{h}])$, the $i_{\mathrm{th}}$ restriction step concatenates $n_d$ convolutions of different dilations,
    \item $\mathcal{T}_i: \mathbb{R}^{c_{\text{in}}\times n_{\text{in}}} \rightarrow \mathbb{R}^{c_{\text{out}}\times n_{\text{out}}}$, the deconvolution (also called transposed convolution) in the $i_{\mathrm{th}}$ prolongation step, with $n_{\text{in}} < n_{\text{out}}$,
    
    \item $\mathcal{P}_i \bm{h}  = \mathcal{N}(\sigma(\mathcal{T}_i \bm{h}))$, the $i_{\mathrm{th}}$ prolongation step,
    
    \item $\mathcal{C}_m: \mathbb{R}^{c\times n} \rightarrow \mathbb{R}^{n}$, a final convolution to mix different channels and produce a one-channel output in the end,
    \item ${\mathcal{V}} = {\mathcal{C}}_\text{m} \circ \mathcal{P}_L\circ \mathcal{P}_{L-1} \circ \cdots \circ \mathcal{P}_1 \circ \mathcal{R}_L \circ \cdots \circ \mathcal{R}_2 \circ \mathcal{R}_1 $, the operation $\mathcal{V}$ is the composition of $L$ restriction and $L$ prolongation steps, followed by $\mathcal{C}_m$,  
\end{itemize}
where  $i=1,2,\dots, L$, $j=1,2,\dots, n_d$, $\bm{h}$ denotes hidden representation, $L$ the total number of downsampling and upsampling steps, $n_d$ the number of different dilation parameters. We use the symbol $[\cdot]$ to indicate concatenation of the tensors along the channel dimension. \par

Let $\mathcal{V}^k$ denote the $k_{\mathrm{th}}$ V-shaped operation, $k=1,2,\dots, n_b$, with $n_b$ being the number of blocks. We explicitly include $\bm{\eta}, \mathbf{f}$ in the input for $\mathcal{V}^k$ by using $\mathbf{X}_k = [\mathbf{u}_k, \bm{\eta}, \mathbf{f}]$ as the tensor concatenating $\mathbf{u}_k$, $\bm{\eta}$, $\mathbf{f}$ in the channel dimension. Define the update $\mathbf{u}_k \mapsto \mathbf{u}_{k+1}$ predicted by $\mathcal{V}^k$:
\begin{equation}
\begin{split}
    \mathbf{u}_{k+1} &= \mathbf{u}_k + \mathcal{V}^k (\mathbf{X}_k), \quad \mathbf{u}_0 = \mathbf{0},\\
    \mathbf{X}_{k+1} &= [\mathbf{u}_{k+1}, \bm{\eta}, \mathbf{f}] \in \mathbb{R}^{3\times n}.
\end{split}
\end{equation}
\par

\begin{figure}
  \centering
  \includegraphics[width=0.95\textwidth]{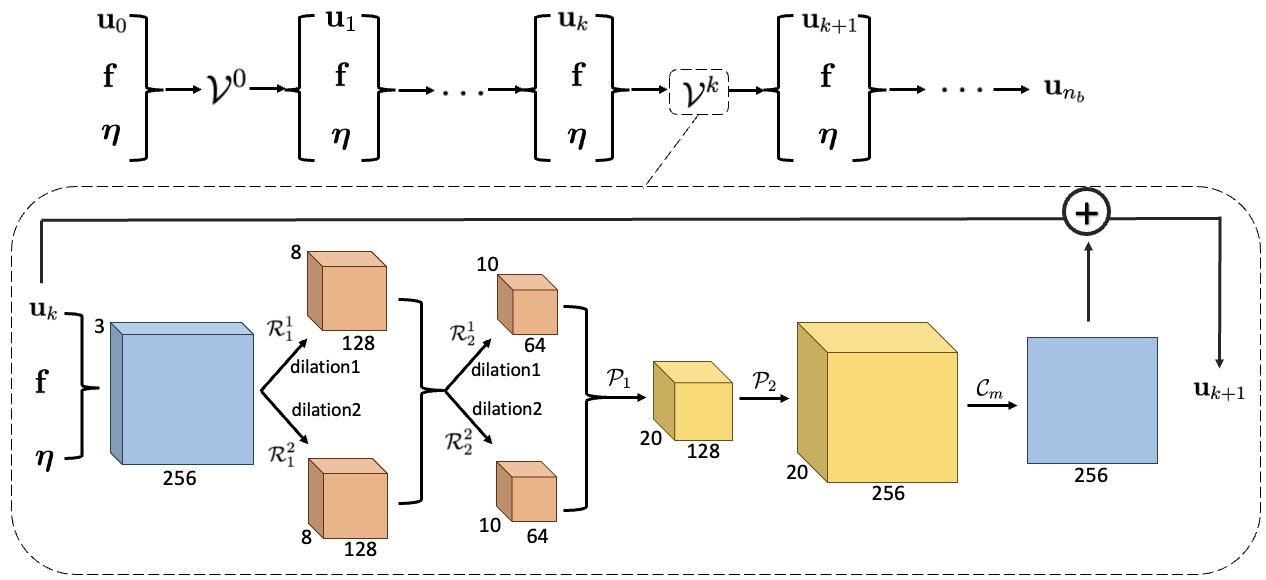} 
  \caption{{The architecture of a block completing one iterative update}. Firstly, concatenate $\mathbf{u}_i, \mathbf{f}, \bm{\eta}$ to form a three-channel input $[\mathbf{u}_i, \mathbf{f}, \bm{\eta}]$, each channel of size $256\times 256$. In the first restriction step, use two paths of convolutions, dilation=1 and dilation=2, denoted by $\mathcal{R}_1^1$ and $\mathcal{R}_1^2$, for downsampling to $128\times128$, and concatenate the output from two paths in the channel dimension. Similarly, in the second restriction step, apply $\mathcal{R}_2^1$ and $\mathcal{R}_2^2$ to arrive at $64\times64$ (the coarsest level).  Then, apply the prolongation operation twice (denoted by $\mathcal{P}_1$, $\mathcal{P}_2$) for upsampling, each time by a factor of 2, to get back to the original mesh size. Finally, a mixing convolution module $\mathcal{C}_m$ fuses the multichannel representation into one channel update. Such an update is added to $\mathbf{u}_i$ to get $\mathbf{u}_{i+1}$. $\mathbf{u}_{i+1}$ is then concatenated with $\mathbf{f}, \bm{\eta}$ again, and the obtained $[\mathbf{u}_{i+1}, \mathbf{f}, \bm{\eta}]$ is fed into the next block. }
  \label{f:fig1}
\end{figure}

In \Cref{f:fig1}, we show an example of one network block with $n_d = 2$ and $L=2$. \xnet\ is created by composing $n_b$ blocks. \par
Green's functions of elliptic operators are non-local. From this perspective, the solution operator should be a non-local operator applied on $\mathbf{f}$. The standard convolution module from computer vision is local. To capture nonlocality we would need several layers within one block, making training harder. \par
Dilated convolutions introduced in  \cite{yuMultiScaleContextAggregation2016} allow for an increased receptive field while using the same number of weights as a narrower convolution kernel. A wider receptive field means that the network can capture more global features from the input field. Directly using larger kernel sizes, however, increases the number of weights in the network, and leads to larger computational overhead and an increased risk of overfitting, especially on small datasets. Dilated convolution enables a wider stencil without increasing the number of weights.\par

While $\bm{\eta}$ and $\mathbf{f}$ are not changing from block to block, including them repeatedly in the input to each block helps with training and generalization empirically.
Since each convolution takes linear work, the complexity of applying \xnet\ would be $O(n_bn)$.
We remark that although the resolution of the prediction $\mathbf{u}$ is fixed, it can be evaluated at any point using some interpolation scheme, such as bilinear interpolation.

\subsection{Loss Function with Spatial Derivatives}
\label{s:dloss}
Many PDE-related downstream applications require accurate prediction of not only  $u$, but also the spatial derivative $\nabla u$. This is also reflected in the fact that $H^1$ norm is required in theoretical error estimates. To ensure low error in $\nabla u$ and enhance continuity in predicting $u$, we add $\nabla u$ to the loss function, guiding the network to also match the spatial derivatives.\par
The loss function is given in \cref{e:loss}, where $s$ is the number of samples, $\mathbf{u}_j^{\ast}$ the $j_{\mathrm{th}}$ true solution, $\mathbf{u}_j$ the network prediction, $\partial_x$ and $\partial_y$ the partial derivatives. Note that since $\mathbf{u}_j^{\ast}, \mathbf{u}_j \in \mathbb{R}^n$, $\|\cdot\|_2$ corresponds to a discretized $L_2$ norm.
\begin{equation}
\label{e:loss}
\mathcal{L}=\frac{1}{s} 
\sum_{j=1}^s \left( \frac{\|\mathbf{u}_{j}-\mathbf{u}_j^{\ast}\|_2}{\|\mathbf{u}_j^{\ast}\|_2} + \frac{\|\partial_x \mathbf{u}_j- \partial_x\mathbf{u}_j^{\ast}\|_2}{2\|\partial_x \mathbf{u}_j^{\ast}\|_2}  + \frac{\|\partial_y \mathbf{u}_j- \partial_y \mathbf{u}_j^{\ast}\|_2}{2\|\partial_y \mathbf{u}_j^{\ast}\|_2}
\right).
\end{equation} 
Empirically, we find it helps to reduce the errors in $\nabla u$ significantly, and also improves accuracy in $u$ slightly, possibly by removing checkerboard artifacts resulted from deconvolution \cite{odenaDeconvolutionCheckerboardArtifacts2016}. We discuss these results in \Cref{s:comp} and \Cref{t:gradu_loss}. 
Relative errors are used for $u$ and $\nabla u$, as we do not want to bias different samples by their norms.

{\bf Training \xnet}:
To be concrete we discuss training for the Poisson problem. 
First we generate $s$ samples $\{\bm{\eta}_i, \bm{f}_i\}$. Then for each sample we solve~\cref{e:gh} to generate $\{\bm{u}_i\}$.
In  the forward  problem  \xnet\ is trained to predict $\bm{u}$ given $\bm{\eta}$, $\bm{f}$.
In the inverse problem  \xnet\ is trained to predict $\bm{\eta}$ given $\bm{u}$, and $\bm{f}$. 
During infrence no additional training takes place. We simple evaluate \xnet\ using the appropriate inputs. 

\section{Related Work}
\label{s:related}
Many alternative methods have been proposed \cite{liFourierNeuralOperator2021, luLearningNonlinearOperators2021, bhattacharyaModelReductionNeural2021, caoChooseTransformerFourier2021, liTransformerPartialDifferential2023, heMgNOEfficientParameterization2023, yePDEformerFoundationModel2024, lanthalerNonlinearReconstructionOperator2022, fanaskovSpectralNeuralOperators2024, linGreenMultigridNetwork2024} for operator learning, and it is beyond the scope of this paper to discuss all of them.
We discuss three related methods: DeepONet \cite{luLearningNonlinearOperators2021} (a network based on separation of variables), and Fourier Neural Operator \cite{liFourierNeuralOperator2021} (a network that builds approximations in spectral space), and U-Net \cite{ronnebergerUNetConvolutionalNetworks2015a} (one of the first multicomponent approximation nets).

\subsection{DeepONet}
{\bf Deep operator net (DeepONet)} \cite{luLearningNonlinearOperators2021} is a deep variant of the operator network proposed in \cite{chenUniversalApproximationNonlinear1995}. Let $G$ be an PDE solution operator taking an input coefficient function $\eta$, with $G(\eta)$ being the corresponding solution function. For any point $y$ in the domain of $G(\eta)$, $G(\eta)(y) \in \mathbb{R}$. Following a separation of variables expansion \cite{lanthalerNonlinearReconstructionOperator2022, lanthalerErrorEstimatesDeepOnets2022, kopanicakovaDeepOnetBasedPreconditioning2024}, the operator is expanded as:

\begin{equation}
\label{e:deeponet}
G(\eta)(y) \approx \sum_{i=1}^p \underbrace{b_i(\eta)}_{\text {branch }} \underbrace{t_i(y)}_{\text {trunk }},
\end{equation}
where $b(\cdot)$ and $t(\cdot)$ refer to the branch and trunk networks, $\eta$ the PDE coefficient and the input of the branch net, $y$ the query point and the input of the trunk net, and $p$ is the hyperparameter controlling the number of terms in the expansion.\par
The branch networks $b_i: \mathbb{R}^n \rightarrow \mathbb{R}$ take $\bm{\eta}$ as input, and extracts hidden representations through a multilayer perceptron (MLP), a convolutional neural network (CNN), or some other encoder network \cite{luLearningNonlinearOperators2021, zhangHybridIterativeNumerical2022, kontolatiLearningNonlinearOperators2024}. The trunk network $t_i: \mathbb{R}^d \rightarrow \mathbb{R}$ transforms a query $\mathbf{y}$ through an MLP. Then, these hidden representations go through element-wise multiplication and summation to output the target function $G(\eta)$ evaluated at the query point. The loss function is the mean squared error between the true value of $G(\eta)(y)$ and the network prediction \cite{luLearningNonlinearOperators2021, zhangHybridIterativeNumerical2022}.
DeepONet has been applied in PDEs from different scientific disciplines, such as astrophysics \cite{maoPPDONetDeepOperator2023}, chemistry \cite{linOperatorLearningPredicting2021}, and fluid dynamics \cite{caiDeepMMnetInferring2021}. Furthermore, efforts have been to combine it with classical iterative methods to accelerate the solving process \cite{zhangBlendingNeuralOperators2024, hu_hybrid_2024} and provide preconditioners \cite{kopanicakovaDeepOnetBasedPreconditioning2024}.\par

\subsection{Fourier neural operator}
{\bf Fourier neural operator (FNO) }\cite{liFourierNeuralOperator2021} is based on parametrizing the Green's function of the PDE operator, and approximating the non-local integration through multiplication in Fourier space. The integration kernel is parametrized as a learnable tensor also in Fourier space. Furthermore, the kernel is truncated at  $k_{\text{max}}$ Fourier modes to ensure efficiency. \par
To learn the mapping from coefficient $\eta$ to solution $u$, FNO first lifts $\bm{\eta}\in \mathbb{R}^n$ to multiple channels by a convolution $C_1: \mathbb{R}^n \rightarrow \mathbb{R}^{c\times n}$ as $v_0 = C_1\eta$. Then it updates hidden representations $v_0 \mapsto v_1 \mapsto \cdots \mapsto v_K$ with $K$ FNO layers, where each $v_j \in \mathbb{R}^{c\times n}$, $j=0,1,\dots,K$. Each FNO layer is composed of a kernel integral and a residual connection:
\begin{equation}
v_{l+1}=\sigma\left(\mathcal{W} v_l+ \mathcal{F}^{-1}\left(R\cdot\mathcal{F}(v_l)\right)\right)\\
\end{equation}
\begin{equation*}
    \left(R \cdot \mathcal{F}\left(v_l\right)\right)_{i, k}=\sum_{j=1}^c R_{j, i, k}\left(\mathcal{F}\left(v_l\right)\right)_{j, k}, \quad \forall k=1, \ldots, k_{\text{max} }, i=1, \ldots, c,
\end{equation*}
where $\mathcal{F}$ and $\mathcal{F}^{-1}$ denote Fourier transform and its inverse, $R\in \mathbb{C}^{k_{\text{max}}\times c\times c}$ is a complex-valued weight tensor truncated by $k_{\text{max}}$ modes, $\mathcal{W}: \mathbb{R}^{c\times n} \rightarrow \mathbb{R}^{c\times n}$ is a linear transformation, $\sigma: \mathbb{R} \rightarrow \mathbb{R}$ is an element-wise nonlinear activation function. Typically, $\mathcal{W}$ is a convolution, and $\sigma$ is the Gaussian Error Linear Unit (GELU) function. Finally, $v_K$ is projected back to the original dimension by another convolution $C_2: \mathbb{R}^{c\times n}\rightarrow \mathbb{R}^n $ as $\hat{u} = C_2v_K$. $\hat{u}$ is the prediction. The relative error $\frac{\|u - \hat{u}\|_2}{\|u\|_2}$ is used as the loss function. FNO has also been found useful in many applications \cite{wenUFNOEnhancedFourier2022, rahmanUshapedNeuralOperators2023, pathakFourCastNetGlobalDatadriven2022, gopakumarFourierNeuralOperator2023}.\par

Although DeepONet and FNO have been successfully applied in various scenarios, the existing literature predominantly focuses on relatively smooth coefficients with small contrast. In the Darcy flow example from \cite{liFourierNeuralOperator2021, wenUFNOEnhancedFourier2022}, the random coefficient field is binary valued with small contrast. For the solution operator examples from \cite{luLearningNonlinearOperators2021} and \cite{luDeepXDEDeepLearning2021a}, the random fields are restricted to smooth components of certain function spaces. 
More detailed discussions about these issues can be found in \cite{lanthalerNonlinearReconstructionOperator2022, fanaskovSpectralNeuralOperators2024}.
\par

\subsection{PDE-Net}
{\bf PDE-Net} \cite{longPDENetLearningPDEs2018, longPDENetLearningPDEs2019} is the work that inspired \xnet. It is designed to accurately predict system dynamics and to uncover the underlying hidden PDE models. We emphasize that it was not proposed for parametric PDEs. Given $\{u_0, \dots, u_t \}$, where $t$ is the time index, PDE-Net tries to extrapolate in time. In PDE-Net, a $\delta t$-block is trained to advance solution $u(t,x)$ to $u(t+\delta t, x)$. 
\begin{equation}
\hat{u}(t+\delta t, x) = D_0 u(t, x) + \delta t\cdot F(x, D_{00}u, D_{01}u, D_{01}u, D_{20}u, \dots),
\end{equation}
where $\hat{u}(t+\delta t, x)$ is the network prediction, $F$ the nonlinear response function to be learned, and $D_{\{\cdot\}}$ the convolution filters approximating different orders of differentiation.
By stacking multiple $\delta t$-blocks, the network is able to predict long time dynamics. Particularly, each $\delta t$-block uses constrained convolutions to learn differential operators, exploiting the relationship between the orders of differential operators and the orders of sum rules of filters \cite{longPDENetLearningPDEs2018}. Therefore, by inspecting the learned convolution filters, one can discover the form of the underlying PDE from data. In another work \cite{khooSolvingParametricPDE2021}, the authors explore the possibility of using a convolutional network to compute low-dimensional quantities of interest dependent on the PDE solution. An iterative update structure is also incorporated to resemble the solving process.\par
\subsection{U-Net}
           {\bf U-Net}, initially proposed for image segmentation \cite{ronnebergerUNetConvolutionalNetworks2015a}, has become a foundational architecture in computer vision tasks that require pixel-wise prediction. Characterized by its encoder-decoder structure with skip connections, U-Net enables the capture of both low-level and high-level features, making it a highly effective architecture. The success of U-Net has inspired various neural architectures across fields that aim to combine multiscale feature extraction with efficient image processing \cite{zhuBayesianDeepConvolutional2018, guoConvolutionalNeuralNetworks2016, sapplDeepLearningPreconditioners2019, raonicCONVOLUTIONALNEURALOPERATORS2023}.
           The authors in \cite{azulayMultigridAugmentedDeepLearning2023}, use a U-Net approximation with in a flexible GMRES solver with  multigrid preconditioner for the  time-harmonic Helmholtz problem.  The U-net operator is used as coarse grid solver within a V-cycle.
           However, this approach requires generalization over both $\eta$ and the residuals. 
           The latter is completely arbitrary, and require very expensive training and is only valid for specific $\eta$.
           So overall the method seems to be more effective when generalizing over $f$ for known $\eta$.
           We discuss a similar approach for \xnet\ in \Cref{s:ablation}.
           
\subsection{Novelty of \xnet}
\label{s:novelty}
\xnet\ borrows ideas from PDE-Net, ResNet \cite{heDeepResidualLearning2016}, and iterative multigrid solvers. \xnet\ constructs a sequence $u_0$, $u_1$, $u_2$, ..., $u_{n_b}$ through operations of $n_b$ blocks to transform an initial guess $u_0=0$ to the approximation of the true solution. The weights of the blocks can be different, or can be the same (see \Cref{s:sharing} for the parameter sharing case).

PDE-Net \cite{khooSolvingParametricPDE2021},  is specifically designed for time-dependent equations, and not for parametric PDEs.
Since PDE-Net must learn trajectory dynamics, it needs trajectory pairs of $\left\{u_j(t_i), u_j(t_{i+k})\right\}$, where $j$ indicates a certain initial condition, each $t_i$ is a starting time, $t_{i+k}$ is the prediction time. The form of response function $F$ is also needed in order to construct the $\delta t$-blocks. The convolutional filters are  constrained to approximate differential operators of various orders, while in our network there are no constraints. 

In the FNO network, the integral kernel is parameterized in Fourier space, and only a maximal number $k_{\text{max}}$ of Fourier modes is used in parametrization to lower the computational complexity. Our network performs (multilevel) convolutions directly in the physical space, and there is no explicit spectral truncation.

Let us also  mention that there are many other specialized nets for particular instances of forward and inverse problems.  The literature is vast, so we will not discuss it here. But representative examples inlude teh following:
The authors in \cite{khoo2019switchnet} discuss the Helmholtz problem in which the scatterer $\eta$ is a set of point scattered and the map is not to the entire scattered field but the field at a small number of receiver locations and generalizes for a small number of point scatterers in a homogeneous medium.

In conclusion \xnet\ (\cref{f:fig1}) has the following features that improve its performance significantly compared to a U-net:
\begin{enumerate*}[label=(\roman*)]
\item repeated blocks with dimension equal to the input;
\item $\bm{f}$ and $\bm{\eta}$ aree repeated at each block;
\item the two-way dilation branch;
  and 
\item use of the $H_1$ norm in the loss function.
\end{enumerate*}
In the appendix we discuss several other variants that did not make a significant difference in the examples we have considered here. 

\section{Results}
\label{s:sec_results}

In this section, we present numerical studies evaluating the performance of \xnet. As discussed in \Cref{s:related},
we compare with three alternative  methods: Fourier Neural Operator \cite{liFourierNeuralOperator2021}, DeepONet \cite{luLearningNonlinearOperators2021}, and U-Net \cite{ronnebergerUNetConvolutionalNetworks2015a}.
Below we summarize the results:
\begin{itemize}
\item \Cref{t:various}: Comparison of \xnet\  with other methods in test errors for the Poisson, Helmholtz, Darcy problems
\item \Cref{t:converge}: Convergence rate of \xnet\ in the Poisson problem, as a function of the number of training PDEs and network size
\item \Cref{f:fig2}: Convergence plot and error distribution for the Poisson problem, comparing with the POD method 
\item \Cref{t:mesh}: Dependence of \xnet\ performance on the underlying PDE mesh discretization size $n$  for the Poisson problem 
\item \Cref{t:UQ}: Convergence for different QoIs and comparison with  POD 
\end{itemize}

\subsection{Experiment Setup}
If not otherwise specified, the problem domain $\Omega=(0,1)^2$ is discretized on a regular grid of $256\times 256$. 500 training samples and 200 testing samples are used for baseline experiments. Each sample corresponds to solving a PDE with a sample of the random coefficients. 
The number of samples used for training  varies when we perform convergence tests.

The architectures used for DeepOnet and FNO are summarized as follows.
For DeepONet we use a CNN as the branch network, and an MLP as the trunk network. For FNO, we adopt the implementation and architecture configuration from \cite{liFourierNeuralOperator2021} for FNO.
\xnet\ uses ReLU as activation, and batch normalization \cite{ioffeBatchNormalizationAccelerating2015} for regularization. Hyperparameter tuning is done for each problem, and the best performance is reported, which results in slight variations in network size. See \ref{a:training} for more details.\par
For the Poisson and Helmholtz problems, we train for about 4500 epochs with an initial learning rate of 0.001 that is later adjusted by a ReduceLROnPlateau learning rate scheduler \cite{lecunDeepLearning2015}. For the Darcy problems, the network was trained for about 1500 epochs. The number of epochs is set empirically so that the loss no longer decreases. We use \texttt{AdamW} as the optimizer. PyTorch is used for all deep learning scripts. For more details  see \ref{a:training}.\par 
Unless otherwise specified, we report relative errors of the form: $\frac{\|\mathbf{u}-\mathbf{\hat{u}}\|_2}{\|\mathbf{u}\|_2}$, where $\|\cdot\|_2$ is $L_2$ norm, $\mathbf{u}$ is the discretized true solution, $\mathbf{\hat{u}}$ is the discretized model prediction.  For some cases, we also report relative errors of $\nabla \mathbf{u}$. 
\subsection{Test Accuracy}
\label{s:comp}
We list test errors of in \Cref{t:various} for Poisson, Helmholtz, Darcy problems and different contrast levels. The contrast of coefficient $(1+\bm{\eta})$ is defined as $\frac{{\rm max}(1+\bm{\eta})}{{\rm min}(1+\bm{\eta})}$, while the contrast of $\mathbf{a}$ is fixed to be $\frac{12}{3}=4$ (See \cref{e:phi}). We report number of trainable weights for all neural networks. We have performed hyperparameter tuning and report the best case for each network. As a result, the networks have slightly different number of weights.\par
With 500 training samples, \xnet\ is the most accurate in Poisson and Darcy problems. When the contrast gets smaller, all methods obtain lower errors, especially the POD and FNO. 
All methods work well for the Darcy flow problem. For Helmholtz problem we get a different story. Despite the relative low wave numbers and the mild contrast, n
all methods result in large errors. With a budget of 1000 samples, the POD method becomes the best in low-contrast Poisson problem, while \xnet\ is  second best.\par
The training errors are much lower than testing errors for FNO in high-contrast Poisson problems, which implies difficulties in training the FNO net. Both \xnet\ and FNO show signs of overfitting in the Helmholtz problem. We attempted several simple approaches to improve these results, for example weight decay (i.e., $\ell_2$ weight regularization), normalization, $L_2$ and $L_1$ regularization for the ouput, and hyperameter tuning. All failed. \par
Here we speculate on the reasons for the unsatisfactory performance of different methods. U-Net \cite{ronnebergerUNetConvolutionalNetworks2015a} was originally designed for image segmentation tasks. To recognize and label semantically meaningful objects of different sizes in an image, it reuses encoded feature maps using  skip connection, and apply convolutions at various spatial scales. However, the PDE task is different because the original image (coefficient field) is not directly similar to the output, and there can be deleterious information loss during spatial compression. Finally, there is no mechanism to account for long-range interactions between input and output. \par
In the implementation of FNO \cite{liFourierNeuralOperator2021}, the integration kernel is assumed to admit a low-frequency truncation in the spectral domain with a maximal number $k_{\text{max}}$ of Fourier modes. While this is essential for efficiency, it may lead to larger errors. \par
As for DeepONet \cite{luLearningNonlinearOperators2021}, using a trunk net to find 1D hidden space embeddings of $\bm{\eta}$ may destroy important spatial details necessary for solution predictions. The use of MLPs in the trunk net and/or branch net may lead to spectral bias in the model \cite{rahamanSpectralBiasNeural2019, basriFrequencyBiasNeural2020,  tancikFourierFeaturesLet2020a, zhangHybridIterativeNumerical2022, linOperatorLearningPredicting2021, kopanicakovaDeepOnetBasedPreconditioning2024}. Some limitations of DeepONet in handling discontinuous solutions are also discussed in \cite{lanthalerNonlinearReconstructionOperator2022}.
\newcolumntype{P}{>{\centering}p{2.9cm}}
\renewcommand\cellalign{{}{P}}

\begin{table}[]
\centering
\caption{Testing errors and number of model weights (M stands for million) for Poisson, Darcy, Helmholtz problems. For the POD method, testing error means prediction error. We use $\|\bm{\eta}\|_{\infty}$ to indicate the contrast level for coefficient $(1+\bm{\eta})$, since $\frac{\text{max}(1+\bm{\eta})}{\text{min}(1+\bm{\eta})} = \frac{\text{max}(1+\bm{\eta})}{1}\approx\|\bm{\eta}\|_{\infty}$. The contrast in Darcy's coefficient $\mathbf{a}$ is fixed to  $\frac{12}{3}$. }
\resizebox{\textwidth}{!}{
\begin{tabular}{|cc|c|c|c|c|c|}
\hline
                               &                                                                                                            & \xnet\                     & POD                       & U-Net                            & FNO                             & DeepONet                        \\ \hline
\multicolumn{1}{|c|}{}         & \cellcolor[HTML]{C0C0C0}\begin{tabular}[c]{@{}c@{}}Poisson,  \\ $\|\bm{\eta}\|_{\infty}\sim$1e6\end{tabular}               & \cellcolor[HTML]{C0C0C0}\makecell{$\nmm2.49\%\pm0.17\%$\\$1.18$M}   & \cellcolor[HTML]{C0C0C0}24.37\%  & \cellcolor[HTML]{C0C0C0}\makecell{$53.69\%\pm0.83\%$\\$7.76$M}  & \cellcolor[HTML]{C0C0C0}\makecell{$45.05\%\pm0.68\%$\\$2.36$M} & \cellcolor[HTML]{C0C0C0}\makecell{$99.99\%\pm0.01\%$\\$1.14$M} \\
\multicolumn{1}{|c|}{}         & \begin{tabular}[c]{@{}c@{}}Poisson,\\ $\|\bm{\eta}\|_{\infty}\sim$1e3\end{tabular}                                         & \makecell{$\nmm1.29\%\pm0.05\%$\\$1.08$M}                           & 1.64\%                           & \makecell{$48.43\%\pm1.32\%$\\$7.76$M}                          & \makecell{$\nmm7.97\%\pm0.43\%$\\$0.53$M}                          & \makecell{$90.63\%\pm3.82\%$\\$1.20$M}                         \\
\multicolumn{1}{|c|}{$s$ = 500}  & \cellcolor[HTML]{C0C0C0}Darcy                                                                              & \cellcolor[HTML]{C0C0C0}\makecell{$\nmm1.23\%\pm0.23\%$\\$0.29$M}   & \cellcolor[HTML]{C0C0C0}2.15\%   & \cellcolor[HTML]{C0C0C0}\makecell{$\nmm7.08\%\pm0.35\%$\\$7.76$M} & \cellcolor[HTML]{C0C0C0}\makecell{$\nmm2.47\%\pm0.15\%$\\$1.33$M}  &\cellcolor[HTML]{C0C0C0}\makecell{$\nmm5.36\%\pm0.26\%$\\$1.61$M}\\
\multicolumn{1}{|c|}{}         & \begin{tabular}[c]{@{}c@{}}Helmholtz, \\ $(2\pi\kappa)^2 = 100$\end{tabular}                           & \makecell{$24.76\%\pm1.17\%$\\$1.18$M}                          & 148.62\%                         & \makecell{$44.95\%\pm0.85\%$\\$7.76$M}                         & \makecell{$11.58\%\pm0.20\%$\\$0.53$M}                      & \makecell{$96.33\%\pm3.17\%$\\$1.20$M}                         \\
\multicolumn{1}{|c|}{}         & \cellcolor[HTML]{C0C0C0}\begin{tabular}[c]{@{}c@{}}Helmholtz,  \\ $(2\pi\kappa)^2=500$\end{tabular} & \cellcolor[HTML]{C0C0C0}\makecell{$83.21\%\pm0.43\%$\\$1.18$M}  & \cellcolor[HTML]{C0C0C0}129.92\% & \cellcolor[HTML]{C0C0C0}\makecell{$86.80\%\pm1.35\%$\\$7.76$M} & \cellcolor[HTML]{C0C0C0}\makecell{$81.68\%\pm0.45\%$\\$0.88$M}   & \cellcolor[HTML]{C0C0C0}\makecell{$99.99\%\pm0.01\%$\\$1.20$M} \\ \hline
\multicolumn{1}{|c|}{}         & \begin{tabular}[c]{@{}c@{}}Poisson,  \\ $\|\bm{\eta}\|_{\infty}\sim$1e6\end{tabular}                                       & \makecell{$\nmm1.24\%\pm0.19\%$\\$3.35$M}                           & 7.61\%                           & \makecell{$49.49\%\pm0.88\%$\\$7.76$M}                          & \makecell{$36.63\%\pm1.11\%$\\$2.51$M}                         & \makecell{$98.94\%\pm0.01\%$\\$1.14$M}                        \\
\multicolumn{1}{|c|}{}         & \cellcolor[HTML]{C0C0C0}\begin{tabular}[c]{@{}c@{}}Poisson,  \\ $\|\bm{\eta}\|_{\infty}\sim$1e3\end{tabular}               & \cellcolor[HTML]{C0C0C0}\makecell{$\nmm0.76\%\pm0.01\%$\\$1.18$M}   & \cellcolor[HTML]{C0C0C0}0.26\%   & \cellcolor[HTML]{C0C0C0}\makecell{$44.67\%\pm0.49\%$\\$7.76$M}  & \cellcolor[HTML]{C0C0C0}\makecell{$\nmm5.81\%\pm0.55\%$\\$0.75$M} & \cellcolor[HTML]{C0C0C0}\makecell{$95.66\%\pm1.02\%$\\$1.18$M} \\
\multicolumn{1}{|c|}{$s$ = 1000} & Darcy                                                                                                      & \makecell{$\nmm0.78\%\pm0.02\%$\\$0.24$M}                          & 1.47\%                           & \makecell{$\nmm6.91\%\pm0.45\%$\\$7.76$M}                          & \makecell{$\nmm1.30\%\pm0.15\%$\\$1.09$M}                          & \makecell{$\nmm5.24\%\pm0.07\%$\\$1.62$M}                          \\
\multicolumn{1}{|c|}{}         & \cellcolor[HTML]{C0C0C0}\begin{tabular}[c]{@{}c@{}}Helmholtz, \\ $(2\pi\kappa)^2=100$\end{tabular}   & \cellcolor[HTML]{C0C0C0}\makecell{$21.02\%\pm0.43\%$\\$1.40$M} & \cellcolor[HTML]{C0C0C0}148.64\% & \cellcolor[HTML]{C0C0C0}\makecell{$43.74\%\pm0.31\%$\\$7.76$M}  & \cellcolor[HTML]{C0C0C0}\makecell{$11.42\%\pm0.26\%$\\$0.53$M} & \cellcolor[HTML]{C0C0C0}\makecell{$97.62\%\pm2.15\%$\\$1.20$M} \\
\multicolumn{1}{|c|}{}         & \begin{tabular}[c]{@{}c@{}}Helmholtz,\\ $(2\pi\kappa)^2=500$\end{tabular}                           & \makecell{$82.76\%\pm0.51\%$\\$1.97$M}                        & 129.93\%                         & \makecell{$87.76\%\pm0.38\%$\\$7.76$M}                         & \makecell{$80.97\%\pm0.42\%$\\$0.88$M}                          & \makecell{$99.99\%\pm0.01\%$\\$1.24$M}                         \\ \hline
\end{tabular}}
\label{t:various}
\end{table}

\subsection{Convergence with respect to Training Size and Model Size}
\label{s:converge}
\Cref{t:converge} shows errors for the high-contrast Poisson problem obtained with \xnet, with different training size $s$ and network size $M$. An expansion factor $w$ on the number of channels of all the convolution / deconvolution modules is applied to conveniently adjust the network size. For example, if $w=2$ then the number of channels doubles for every convolution operation in \xnet. \par
With only 500 samples, the network can predict solutions with about 2\% error. However, we find that the errors decrease slowly with the size of training data.
Of course due to the statistical sampling of the loss function the expected convergence at best would be $\mathcal{O}(1/\sqrt{s})$.
Since such convergence would require also increasing the network size, the lowest errors are expected to be found along the diagonal of the table, where the increasing dataset size matches network size.

In \xnet, there are two primary ways to introduce more weights: increasing the number of blocks ($n_b$) or expanding the number of channels ($w$) in the convolutional filters, corresponding to increasing the network's depth or width, respectively. We found that, after the number of blocks exceeds a certain value, increasing the width is preferable, as it effectively enhances the model's capability without significantly increasing the training burden.\par
\begin{figure}
  \centering
  \includegraphics[width=0.9\textwidth]{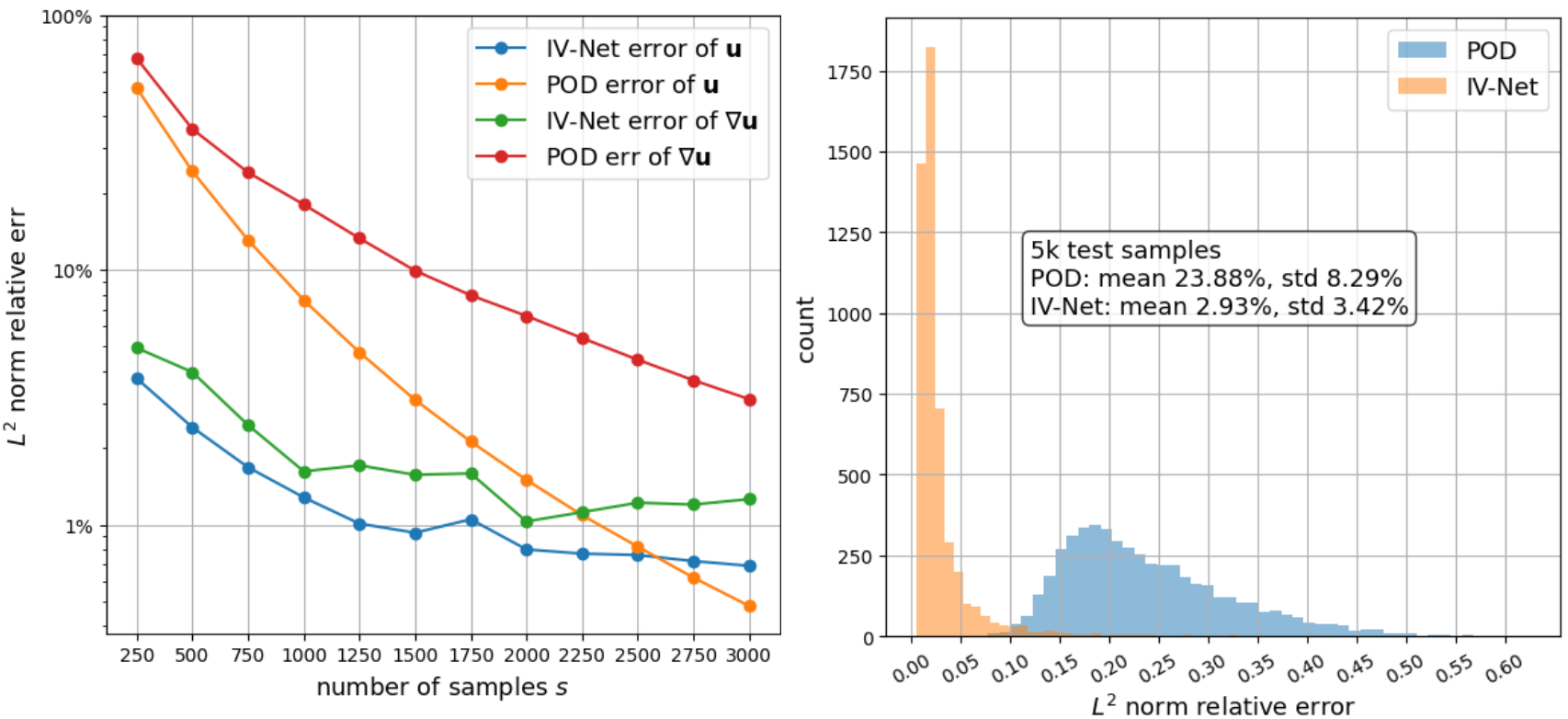} 
  \caption{{Test errors for the high contrast Poisson problem ($\|\bm{\eta}\|_{\infty} \sim $1e6).} Left: error convergence behavior of our network and the POD method with increasing training size. Right: error distributions of \xnet\ and the POD method.}
  \label{f:fig2}
\end{figure}

\Cref{f:fig2}  further contrasts the convergence behavior of our network with that of the POD method. While our network significantly outperforms the projection method when the sample size is small, the network's error saturates as the sample size increases. In contrast, the error of the POD method continues to decrease, with a crossover point occurring at $s = 2500$. For \xnet\, reducing the  error below 1\% becomes quite expensive.
The convergence of the POD method is much more favorable.

\begin{table}[]
\centering
\caption{Convergence of \xnet\ test errors for the Poisson problem with $\|\bm{\eta}\|_{\infty} \sim $1e6. We vary the number of training samples $s$ and network size $M$. $w$ is a multiplicative factor which controls the width of the network. In all results, we use 12 blocks.}
\label{t:converge}
\begin{tabular}{|c|c|c|c|c|}
\hline
\diagbox[]{$s$}{$M$} & \begin{tabular}[c]{@{}l@{}}1.18M \\ ($w$=1.0)\end{tabular} & \begin{tabular}[c]{@{}l@{}}1.97M \\ ($w$=1.3)\end{tabular} & \begin{tabular}[c]{@{}l@{}}2.63M \\ ($w$=1.5)\end{tabular} & \begin{tabular}[c]{@{}l@{}}3.35M\\ ($w$=1.7)\end{tabular} \\ \hline500        & 2.41\%                                                   & 2.54\%                                                   & 2.61\%                                                   & 2.44\%                                                  \\ 
\cellcolor[HTML]{C0C0C0}1000       & \cellcolor[HTML]{C0C0C0}1.75\%                                                   &\cellcolor[HTML]{C0C0C0}1.57\%                                                   & \cellcolor[HTML]{C0C0C0}1.40\%                                                   & \cellcolor[HTML]{C0C0C0}1.28\%                                                  \\ 
2000       & 1.15\%                                                   & 1.00\%                                                   & 0.80\%                                                   & 0.93\%                                                  \\
\cellcolor[HTML]{C0C0C0}3500       &\cellcolor[HTML]{C0C0C0}0.85\%                                                   &\cellcolor[HTML]{C0C0C0}0.72\%                                                   &\cellcolor[HTML]{C0C0C0}0.78\%                                                   &\cellcolor[HTML]{C0C0C0}0.65\%                                                  \\ 
5000       & 0.84\%                                                   & 0.64\%                                                   & 0.50\%                                                   & 0.49\%                                                  \\ \hline
\end{tabular}
\end{table}

\subsection{Convergence Behavior for Different PDE Discretization Sizes (spatial resolution)}
\label{s:mesh}

Using one's favorite inteprolation scheme for inputs and outputs, \xnet\ can be trained in arbitrary resolution related to the underlying PDE.
Furthermore using postprocessing, again  using interpolation, it can be evaluated to arbitrary resolution.
Here, as in the entire paper, we're interested in the training behavior as we fix the \xnet's resolution to be identical to the underlying PDE resolution.

\Cref{t:mesh} presents test errors and the network configurations for various combinations of training set size and spatial resolution. For the coarse spatial resolutions, we find that using two downsampling levels ($L=2$) in \xnet\ leads to overfitting. Therefore, only one level ($L=2$) is used. With this modification, consistently low errors are achieved across different resolutions.\par


\begin{table}[]
\centering
\caption{\xnet\ configuration for different spatial resolutions. Tuples indicate network hyperparameters in the order of $(L,n_b,w)$, followed by relative errors. One downsampling level ($L=1$) means that the images are downsampled by a factor of 2 only once. }
\label{t:mesh}
\begin{tabular}{|c|cc|cc|cc|}
\hline
$n$        & \multicolumn{2}{c|}{$s$ = 500}                   & \multicolumn{2}{c|}{$s$ = 1000}          & \multicolumn{2}{c|}{$s$ = 2000}          \\
\hline
$64\times 64$   & \multicolumn{1}{c|}{$(1,22,0.8)$} & 3.69\%         & \multicolumn{1}{c|}{$(1,14,0.8)$} & 2.40\% & \multicolumn{1}{c|}{$(1,18,0.8)$} & 3.08\% \\ 
\cellcolor[HTML]{C0C0C0}$128\times 128$ & \multicolumn{1}{c|}{\cellcolor[HTML]{C0C0C0}$(1,18,1.0)$} & \cellcolor[HTML]{C0C0C0}5.08\% & \multicolumn{1}{c|}{\cellcolor[HTML]{C0C0C0}$(1,22,1.0)$} & \cellcolor[HTML]{C0C0C0}2.85\% & \multicolumn{1}{c|}{\cellcolor[HTML]{C0C0C0}$(1,22,1.0)$} &\cellcolor[HTML]{C0C0C0}1.92\% \\ 
$256\times 256$ & \multicolumn{1}{c|}{$(2,12,1.0)$}     & 2.41\%         & \multicolumn{1}{c|}{$(2,12,1.7)$}     & 1.28\% & \multicolumn{1}{c|}{$(2,12,1.5)$}     & 0.80\% \\ 
\cellcolor[HTML]{C0C0C0}$512\times 512$ & \multicolumn{1}{c|}{\cellcolor[HTML]{C0C0C0}$(2,12,1.0)$}     & \cellcolor[HTML]{C0C0C0}2.21\%        & \multicolumn{1}{c|}{\cellcolor[HTML]{C0C0C0}$(2,12,1.3)$}     & \cellcolor[HTML]{C0C0C0}1.37\% & \multicolumn{1}{c|}{\cellcolor[HTML]{C0C0C0}$(2,12,1.7)$}     &\cellcolor[HTML]{C0C0C0}1.62\% \\

\hline
\end{tabular}
\end{table}


\subsection{Performance of \xnet\ in Uncertainty Quantification Tasks}
\label{s:UQ}
So far, we have discussed the accuracy of \xnet\ for PDE solves. In this section, we discuss its accuracy in the context of uncertainty quantification of some example quantities of interest (QoIs). \par
Here, we study three QoIs directly related to the solution of the PDE, which are $\|\mathbf{u}\|_{2}$, $\|\mathbf{u}\|_{1}$, and $\sum_{\Omega}\left|\partial_x \mathbf{u}\right|$. We demonstrate the accuracy of \xnet\ in predicting the mean ($\mu$), variance ($\sigma$), and probability of extreme events ($p_{\text{ext}}$) for these QoIs. We first randomly sample $\{\xi_i\}_{i=1}^{s}$, compute the solutions $\{\mathbf{u}_i\}_{i=1}^s$ for the PDEs, and then for each QoI $q$, compute the statistics by:
\begin{itemize}
    \item $\mu(q) = \frac{1}{s}\sum_{i=1}^{s}q(\mathbf{u}_i)$,
    \item $\sigma(q) = \sqrt{\frac{\sum_{i=1}^s \left(q(\mathbf{u}_i) - \mu(q)\right)^2}{s-1}}$,
    \item $p_{\text{ext}}(q) = \frac{1}{s}\sum_{i=1}^s \mathbf{1}\left(|q(\mathbf{u}_i) - \mu(q)| > 3\sigma(q)\right)$,
\end{itemize}
where $\mathbf{1}(t)$ equals 1 if $t$ is true, and 0 otherwise. For comparison, we also list the estimations made by the POD method.\par

\begin{table}
\caption{Comparison of our network and the POD method in predicting statistics of three QoIs. $p_{\mathrm{ext}}$ refers to probability of extreme events where the quantity is more than $3\sigma$ away from mean value. Relative errors are used to present performance results. Positive values indicate overestimation, and negative values underestimation. }
\label{t:UQ}
\centering
\resizebox{\textwidth}{!}{
\begin{tabular}{|c|c|c|c|c|c|c|c|}
\hline \multicolumn{2}{|c|}{} & \begin{tabular}{c} 
Exact\\ 
($s=16000$) \end{tabular}  & \begin{tabular}{c} 
\xnet\\ 
($s=250$)
\end{tabular} & \begin{tabular}{c} 
POD \\
($s=250$)
\end{tabular} & \begin{tabular}{c} 
POD \\
($s=100$)
\end{tabular} & \begin{tabular}{c}
\xnet\\
($s = 500$)
\end{tabular} & \begin{tabular}{c} 
POD \\
($s=500$)
\end{tabular} \\
\hline \multirow{3}{*}{$\|\mathbf{u}\|_{2}$} &\cellcolor[HTML]{C0C0C0}$\mu$ & \cellcolor[HTML]{C0C0C0}$4.98 \times 10^4$ & \cellcolor[HTML]{C0C0C0}$\nm1.81 \%$ & \cellcolor[HTML]{C0C0C0}$-40.76 \%$ & \cellcolor[HTML]{C0C0C0}$-77.91 \%$ & \cellcolor[HTML]{C0C0C0}$-0.80 \%$ & \cellcolor[HTML]{C0C0C0}$-15.86 \%$ \\
& $\sigma$ & $3.38\times 10^4$ & $\nm47.04 \%$ & $-25.74 \%$ & $-64.50 \%$ & $-6.21 \%$ & $-9.76 \%$ \\
&\cellcolor[HTML]{C0C0C0}\begin{tabular}{c} 
$p_{\text{ext}}$
\end{tabular} & \cellcolor[HTML]{C0C0C0}\; $1.82 \times 10^{-2}$ & \cellcolor[HTML]{C0C0C0}$-62.09 \%$ & \cellcolor[HTML]{C0C0C0}$\nm5.49 \%$ & \cellcolor[HTML]{C0C0C0}$\nm18.13 \%$ & \cellcolor[HTML]{C0C0C0}$\nm3.30 \%$ & \cellcolor[HTML]{C0C0C0}$\nm1.65 \%$ \\
\hline 
\multirow{3}{*}{$\|\mathbf{u}\|_{1}$} & $\mu$ & $4.46\times 10^6$ & $\nm2.69 \%$ & $-43.72 \%$ & $-77.13 \%$ & $-0.67 \%$ & $-19.06 \%$ \\
& \cellcolor[HTML]{C0C0C0}$\sigma$ & \cellcolor[HTML]{C0C0C0}$3.55 \times10^6$ & \cellcolor[HTML]{C0C0C0}$\nm47.04 \%$ & \cellcolor[HTML]{C0C0C0}$-34.37 \%$ & \cellcolor[HTML]{C0C0C0}$-69.86 \%$ & \cellcolor[HTML]{C0C0C0}$-9.30 \%$ & \cellcolor[HTML]{C0C0C0}$-15.21 \%$ \\
& \begin{tabular}{c} 
$p_{\text{ext}}$
\end{tabular} &\; $1.63 \times 10^{-2}$ & $-57.05 \%$ & $\nm4.91 \%$ & $\nm30.67 \%$ & $\nm9.81 \%$ & $\nm3.68 \%$ \\
\hline 
\multirow{3}*{$\sum_{\Omega}\left|\partial_x \mathbf{u}\right|$} & \cellcolor[HTML]{C0C0C0}$\mu$ & \cellcolor[HTML]{C0C0C0}$2.84 \times10^5$ & \cellcolor[HTML]{C0C0C0}$\nm1.76 \%$ & \cellcolor[HTML]{C0C0C0}$-36.62 \%$ & \cellcolor[HTML]{C0C0C0}$-74.40 \%$ & \cellcolor[HTML]{C0C0C0}$\nm0.35 \%$ & \cellcolor[HTML]{C0C0C0}$-11.62 \%$ \\
& $\sigma$ & $1.34 \times10^5$ & $\nm35.07 \%$ & $-11.94 \%$ & $-51.57 \%$ & $\nm0.01 \%$ & $-0.75 \%$ \\
& \cellcolor[HTML]{C0C0C0}\begin{tabular}{c} 
$p_{\text{ext}}$
\end{tabular} & \cellcolor[HTML]{C0C0C0}\; $1.59\times 10^{-2}$ & \cellcolor[HTML]{C0C0C0}$-54.09 \%$ & \cellcolor[HTML]{C0C0C0}$\nm 15.72 \%$ & \cellcolor[HTML]{C0C0C0}$\nm33.96 \%$ & \cellcolor[HTML]{C0C0C0}$-1.89 \%$ & \cellcolor[HTML]{C0C0C0}$-0.63 \%$ \\
\hline
\end{tabular}}
\end{table}

From \Cref{t:UQ}, it is evident that our network and the POD method produce close estimates of the statistics of three QoIs with only 500 training samples. \xnet\ is more accurate in estimating mean values and variances, while the POD method performs better in predicting extreme-event probabilities. The POD seems to predict a shifted distribution, as it underestimates the mean but captures the right portion of tail distributions. When the available data size is halved ($s = 250$), even though \xnet\ does well in predicting mean values, it severely overestimates variances and therefore underestimates extreme-event probabilities.
\par
It is important to note that every prediction by the POD method requires solving a reduced linear system by iterative schemes, which can be costly. To match the computational cost, results of the POD method with $s=100$ are also listed. 

\section{Alternative formulations and extensions}
\label{s:sec_diss}
In the following, we discuss whether alternative networks could be used for predicting low-dimensional quantities of interest in \Cref{s:direct_qoi}. As a further application of \xnet, the inverse problem of predicting coefficients from solutions is explored in \Cref{s:inv}. Additionally, a detailed analysis of the \xnet\ architecture through ablation studies is presented in \Cref{s:ablation}. Lastly, we discuss experiment results using a special variant of \xnet\ with identical blocks in \Cref{s:sharing}.
\subsection{Predicting QoI Directly}
\label{s:direct_qoi}
Predicting the quantity of interest $\mathbf{q}\in \mathbb{R}^{d_q}$ seems like an easier task than predicting $\mathbf{u}\in \mathbb{R}^n$ everywhere on a $\sqrt{n} \times \sqrt{n}$ grid since $d_q \ll n$. 
Many neural architectures can serve as a feature extractor to learn a representation that can then be passed to an MLP to predict $\mathbf{q}$. The question is, if we are only interested in $\mathbf{q}$, can we train a simpler network to directly predict $\mathbf{q}$? We investigate this question empirically for the high contrast Poisson problem.\par
To generate a QoI, we fix a random Gaussian matrix $\mathbf{Q}\in \mathbb{R}^{2\times n}$ such that $\mathbf{q} = \mathbf{Qu}$, $\mathbf{q}\in \mathbb{R}^2$. Then, different architectures are tested to predict $\mathbf{q}$ either by $\bm{\xi} \rightarrow \mathbf{q}$ or $\bm{\eta} \rightarrow \mathbf{q}$, where $\bm{\xi}$ is the intrinsic representation of $\bm{\eta}$. 

For mapping $\bm{\xi}$ to $\mathbf{q}$, we tried both MLP and a graph convolution network. The latter treats ring centers as graph nodes, and use coordinates to generate embeddings for the nodes. For mapping $\bm{\eta}$ to $\mathbf{q}$, we tried vanilla CNN, ResNet \cite{heDeepResidualLearning2016}, as well as a Vision Transformer \cite{dosovitskiyImageWorth16x162021}. Regardless of the model adjustments or the amount of training data used (up to 10k), all tested networks show substantial overfitting, with test errors exceeding 50\%. This suggests that these architectures failed to learn effective representations for predicting $\mathbf{q}$, highlighting the challenge of directly predicting low-dimensional QoIs in the high-contrast Poisson problem.

\subsection{Using \xnet\ for Inverse Problem}
\label{s:inv}
To further demonstrate \xnet's potential in PDE-related dense regression tasks, we evaluate its performance in an inverse setting: predicting the coefficient field $\bm{\eta}$ from the solution $\mathbf{u}$. In practical scenarios, a full observation of the accurate $\mathbf{u}$ is not typical, but this setup shows that the model can be easily applied to the inverse problem. \par
Prediction errors of $\bm{\eta}$ are shown in \Cref{t:inv} for the Poisson and Helmholtz problems. Even though the $\bm{\eta}$ fields are complicated, simply swapping $\mathbf{u}$ and $\bm{\eta}$ in the datasets gives good results, without modifying the model architecture. This further demonstrates the flexibility of \xnet\ in handling different tasks related to PDEs. From \Cref{t:inv} we also notice the accuracy of the predictions improves significantly when dealing with coefficients of smaller contrasts. This likely explains why the Helmholtz inverse problem yields lower errors compared to Poisson problems with larger contrast. \par

\begin{table}[]
\centering
\caption{IV-Net for the inverse problem: mapping $\bm{u}$ to $\bm{\eta}$. The two errors in each entry are training and testing relative errors of $\bm{\eta}$.}
\label{t:inv}
\begin{tabular}{|c|c|c|c|}
\hline \xnet\ & \begin{tabular}{c} 
Poisson \\
$\|\bm{\eta}\|_{\infty} \sim$1e6
\end{tabular} & \begin{tabular}{c} 
Poisson \\
$\|\bm{\eta}\|_{\infty} \sim$1e3
\end{tabular} & \begin{tabular}{c} 
Helmholtz \\
$(2\pi\kappa)^2=100$
\end{tabular} \\
\hline $s =5 0 0$ & $2.78 \%, 15.00 \%$ & $1.70 \%, 2.37 \%$ & $0.30 \%, 0.65 \%$ \\

$s =1 0 0 0$ & $2.63 \%, 11.59 \%$ & $1.04 \%, 1.55 \%$ & $0.33 \%, 0.39 \%$ \\
\hline
\end{tabular}
\end{table}

\subsection{Ablation Studies}
\label{s:ablation}
\begin{figure}
  \centering
  \includegraphics[width=0.9\textwidth]{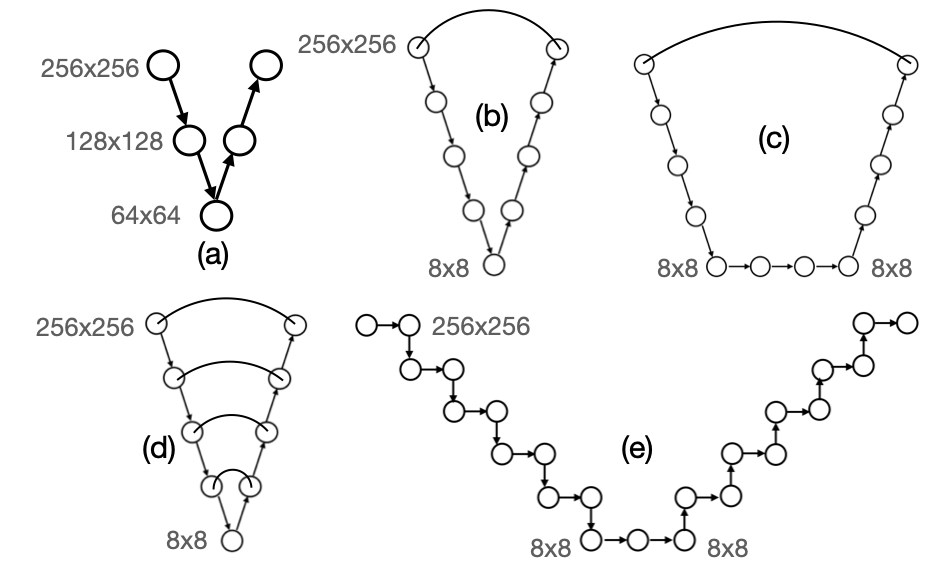} 
  \caption{{Explored structural variants of a single block.}(a) No residual connection; (b) More levels; (c) Coarsest level convolutions; (d) More residual connections; (e) Convolutions for smoothing. Each circle represents an intermediate state. Arcs connecting two states indicate residual connections. Arrows indicate convolutional operations: horizontal ones are unstrided, other ones are strided for downsampling/upsampling purposes. The mesh resolutions are shown next to the states.}
  \label{f:var}
\end{figure}

It is well known that multigrid methods are related to convolutional neural networks \cite{heMgNetUnifiedFramework2019, chenMetaMgNetMetaMultigrid2020, sapplDeepLearningPreconditioners2019}. In \xnet\ , although we do not compute the residual $\mathbf{f}-\mathbf{A}\mathbf{u}_i$ or form traditional prolongation or restriction operators explicitly like \cite{heMgNOEfficientParameterization2023, heMgNetUnifiedFramework2019}, each block in our network does have a multilevel structure that leverages information from different resolutions. It is reasonable to ask if a closer alignment  with traditional multigrid method would help improve performance. In this spirit, we present ablation studies that explore the following architectural variants, many of which drawing ideas from multigrid methods:
\begin{itemize}
\item No residual connection (\Cref{f:var}(a)): remove residual connections;
\item More levels (\Cref{f:var}(b)): use a deeper hierarchy of multi-resolution grids;
\item Coarsest level convolutions (\Cref{f:var}(c)): at the coarsest level, add multiple unstrided convolutions (with activations) to mimic the coarsest level solve of a multigrid solver;
\item More residual connections (\Cref{f:var}(d)): adding residual connections at each level;
\item Convolutions for smoothing (\Cref{f:var}(e)): instead of using strided convolution for direct downsampling or upsampling, add unstrided convolutions before or after that to resemble the pre- and post-smoothing steps;
\end{itemize}

We find that the error increases if we remove residual connections; adding more does not significantly reduce the error. With more levels, the network is more likely to overfit, as is already known from our resolution dependence results (see \Cref{s:mesh}). Finally, adding more convolutions for coarse solve or smoothing does not make much difference. These findings suggest that the current architecture is expressive enough for the problems we considered, while excessive complexity or modification does not necessarily translate to improved results.

\subsection{Identical Blocks}
\label{s:sharing}
One noteworthy variant, or a special case of our network, is to use identical blocks, which means stacking one block $n_b$ times. In \Cref{t:id} we show test errors on coercive problems (Poisson and Darcy) using such a network with $n_b=16$. Moreover, the number of channels are also reduced, resulting in a network of only 14k parameters. Despite the significantly reduced size, the network still produces small errors, demonstrating its effectiveness in a highly parameter-efficient form.

\section{Conclusion}
\label{s:conclusion}
\xnet\ approximates the solution iteratively, computing the update at each step through a V-shaped convolutional block. As an end-to-end dense regression model, \xnet\ exhibits strong performance on high-contrast, high-frequency elliptic equations, achieving competitive predictive accuracy and uncertainty quantification with relatively limited training data. The principal limitations of this study are summarized below.
\begin{itemize}
    \item \xnet fails to reliably solve higher-frequency Helmholtz problems, although the performance on the corresponding inverse problems is still reasonably acceptable.
    \item We only consider problems posed on two-dimensional regular grids. In three dimensions, we anticipate qualitatively similar performance at increased computational cost. The same conceptual framework could, in principle, be extended to unstructured meshes via graph neural networks; however, the detailed development of such an extension remains an open research problem.
    \item The error convergence rate is relatively slow, which renders the method computationally expensive and of limited practicality when higher-order accuracy is required.
    \item Standard convolutional layers, such as those employed in \xnet, may introduce aliasing artifacts \cite{raonicCONVOLUTIONALNEURALOPERATORS2023, karrasAliasFreeGenerativeAdversarial2021}. Mitigating these aliasing effects could further enhance predictive quality.
    \item We compare \xnet\ against a fast GPU-accelerated $O(n^2)$ PDE solver; however, the effectiveness and scalability of the approach in fully three-dimensional settings remain to be systematically evaluated.
    \item \xnet\ can be straightforwardly modified to incorporate FNO blocks or pure filtering operators, but this direction has not been pursued in the present work.
    \item We observe that using the $H_1$ norm in the loss function facilitates training. An even more effective strategy, particularly for inverse problems, is to include derivatives with respect to the input fields in the loss, as in DINO nets~\cite{oleary-roseberryDerivativeInformedNeuralOperator2024}. This derivative-informed approach has not been explored here.
\end{itemize}

We conducted extensive numerical experiments to evaluate the convergence rate of the networks, their dependence on the underlying discretization mesh, and their sensitivity to hyperparameter choices. An important finding is that for comparatively simple benchmark problems, such as the Darcy problem that is widely studied in the literature, none of the considered neural operator architectures consistently outperforms a reference POD (Proper Orthogonal Decomposition) scheme. The POD method is not only the only approach in our comparison that provides practical error guarantees, but it is also essentially straightforward to implement. 

Another key observation is that none of the networks exhibits satisfactory performance for the forward Helmholtz problem. Moreover, although aggregate error statistics (e.g., mean and median errors) appear stable across runs, the maximum error observed at inference time can be substantially large. In addition, in our settings, all networks struggled to reduce the relative error below 1\%. Finally, focusing on lower-dimensional quantities of interest, instead of the full-field solution, tends to make the training procedure more challenging.
Overall, our results suggest that, despite recent advances in neural operators, achieving robust and reliable performance—even for linear two-dimensional problems—while simultaneously surpassing the efficiency and practicality of established reduced-order modeling techniques remains a significant challenge.

\appendix

\section{Data Preparation}
All PDEs are discretized with a second-order finite difference method. For Poisson problems, the right-hand side $f$ and boundary term $g$ are fixed and generated by inserting $u = 1.4\text{cos}(5\pi x)\text{sin}(0.8\pi y).$ For Helmholtz problems, the source term $f$ is fixed to be a sum of three Gaussians with centers $(0.5,0.8), (0.1,0.3), (0.7,0.4)$ and $\sigma = \frac{0.2}{\sqrt{2}},\frac{0.3}{\sqrt{2}},\frac{0.3}{\sqrt{2}}$, and the boundary term $g=0$.
\section{Training Details}
\label{a:training}
The U-Net we used is a standard one, which first downsamples the coefficient from resolution of $256\times 256$ to $16\times 16$ in 4 steps, each step by a factor of 2, and then uses transposed convolutions for upsampling. At each level of resolution (except the coarsest one), a skip connection path concatenates features in the encoder to features in the decoder of corresponding resolution.\par
For training FNO, we refer to the literature \cite{liFourierNeuralOperator2021} for suggested hyperparameter combinations, and conduct a grid search around those values. Particularly, we tried
$\text{num\_modes} = \left\{12, 14, 16, 18, 20, 22, 24, 30\right\}$, $\text{hidden\_channels}= \left\{16, 20, 24, 32, 48 \right\}$, $\text{lifting\_channels}= \left\{16, 32, 64 \right\}$, $\text{projection\_channels}= \left\{16, 32, 64 \right\}$. We did not experiment with larger parameters, as severe overfitting was already observed; increasing the model size appeared to further degrade performance. The input coefficient is concatenated with grid coordinates to ease learning \cite{liFourierNeuralOperator2021, kovachkiNeuralOperatorLearning}.
For training DeepONet, we use the mean squared error loss function used in \cite{luLearningNonlinearOperators2021, zhangHybridIterativeNumerical2022, maoPPDONetDeepOperator2023, caiDeepMMnetInferring2021}. We tried three different ways to treat coefficient functions $\eta$ in the branch net: (1) sample a number of points in the domain and use MLP; (2) treat $\eta$ as images and use CNN; (3) use intrinsic variables $\xi$ (which generate $\eta$) and use MLP. The second approach is the best in our experiments, and this choice aligns with \cite{luComprehensiveFairComparison2022, zhangHybridIterativeNumerical2022}. We tried $\text{hidden\_dim}=\left\{80, 120, 160, 200, 300, 350, 400, 450, 500\right\}$ for the trunk net.\par 
All computations are carried out either on an Nvidia A100 GPU with 40GB memory, or on an Nvidia H100 GPU with 80GB memory.


\section{Errors with Identical Blocks}
In \Cref{t:id}, we report  errors using identical blocks in \xnet.
\begin{table}[]
\centering
\caption{Errors of non-oscillatory problems obtained by using identical blocks. The network has only 14k parameters.}
\label{t:id}
\begin{tabular}{|c|c|c|}
\hline
                        & $s=500$ & $s=1000$ \\ \hline
\begin{tabular}[c]{@{}c@{}}Poisson,  $\|\bm{\eta}\|_{\infty} \sim$1e6\end{tabular} & 10.02\%            & 6.03\%              \\ \hline
\begin{tabular}[c]{@{}c@{}}Poisson,  $\|\bm{\eta}\|_{\infty}\sim$1e3\end{tabular} & 6.44\%             & 4.37\%              \\ \hline
Darcy                                                               & 3.01\%             & 2.35\%              \\ \hline
\end{tabular}
\end{table}

\begin{table} 
\centering
\caption{Errors of the high contrast Poisson problem, using a loss function with and without derivative information. Here \emph{error of $\nabla\mathbf{u}$} means the average of errors on $\nabla_x \mathbf{u}$ and $\nabla_y \mathbf{u}$.}
\begin{tabular}{|cc|c|c|}
\hline
\multicolumn{2}{|c|}{}                                            & with $\nabla \mathbf{u}$   & without $\nabla \mathbf{u}$ \\ \hline
\multicolumn{1}{|c|}{\multirow{2}{*}{s = 500}}  & error of $\mathbf{u}$      & 2.41\% & 3.22\%  \\ \cline{2-4} 
\multicolumn{1}{|c|}{}                          & error of $\nabla \mathbf{u}$ & 3.97\% & 11.45\% \\ \hline
\multicolumn{1}{|c|}{\multirow{2}{*}{s = 1000}} & error of $\mathbf{u}$      & 1.28\% & 2.17\%  \\ \cline{2-4} 
\multicolumn{1}{|c|}{}                          & error of $\nabla \mathbf{u}$ & 1.62\% & 7.43\%  \\ \hline
\end{tabular}
\label{t:gradu_loss}
\end{table}

\section{Adding Spatial Derivatives to Loss}
We show in \Cref{t:gradu_loss} the effect of including relative errors of the spatial derivatives $\nabla \mathbf{u}$ in the loss function \cref{e:loss}. The loss function is also replicated here:
\begin{equation*}
\mathcal{L}=\frac{1}{s} 
\sum_{j=1}^s \left( \frac{\|\mathbf{u}_{j}-\mathbf{u}_j^{\ast}\|_2}{\|\mathbf{u}_j^{\ast}\|_2} + \frac{\|\partial_x \mathbf{u}_j- \partial_x\mathbf{u}_j^{\ast}\|_2}{2\|\partial_x \mathbf{u}_j^{\ast}\|_2}  + \frac{\|\partial_y \mathbf{u}_j- \partial_y \mathbf{u}_j^{\ast}\|_2}{2\|\partial_y \mathbf{u}_j^{\ast}\|_2}
\right).
\end{equation*} 

\section{Scaling \xnet}
We have seen that the sharp contrast in $\eta$ significantly adds to the difficulty of solving the Poisson problem. Could this be adjusted during data preparation to make learning less challenging? \par
The original equation is $\mathbf{Au} = \mathbf{f}$ where $\mathbf{A}$ has $\bm{\eta}$ field in the diagonal. Let $\bm{D} = {\rm diag}((1+\bm{\eta}) - \bm{\Delta})$ be a diagonal matrix, and rewrite the equation as 
\begin{equation}
\bm{D}^{-\frac{1}{2}}\mathbf{A}\bm{D}^{-\frac{1}{2}}\bm{D}^{\frac{1}{2}}\mathbf{u} = \bm{D}^{-\frac{1}{2}}\mathbf{f}.    
\end{equation}
We now construct the dataset with 
\begin{equation}
(\Tilde{\bm{\eta}}, \Tilde{\mathbf{u}}, \Tilde{\mathbf{f}}) = (\bm{D}^{-\frac{1}{2}}\bm{\eta} \bm{D}^{-\frac{1}{2}}, \bm{D}^{\frac{1}{2}}\mathbf{u}, \bm{D}^{-\frac{1}{2}}\mathbf{f}),    
\end{equation}
where $\|\Tilde{\bm{\eta}}\|_{\infty}$ is greatly reduced than $\|{\bm{\eta}}\|_{\infty}$. However, in our experiment, the network produces similar errors after such scaling.

\section{Combining \xnet\ with Krylov Methods}
\label{s:krylov}
As it was done in \cite{azulayMultigridAugmentedDeepLearning2023}, the idea is to try to use a neural operator architecture with a traditional solver to improve convergence.
Define right preconditioner $P$ for solving $Ax=f$: $PAx = Pf$ where $P \approx A^{-1}$. This requires $P$ to be linear and symmetric. Furthermore, applying the preconditioner to a vector $r$ is equivalent to solving $Ax = r$, where $r$ is the residual of the current iteration.
If we try to use \xnet\ to learn the operator $P_{IV}[\eta]$, leading to the following observations:
\begin{enumerate}[nosep]
    \item[{(1)}]  $P_{IV}[\eta]$ is nonlinear and its linearization is non-symmetric,
    \item[{(2)}]  $P_{IV}[\eta]$ is trained on the distribution of $f$, not the distribution of $r$.
\end{enumerate}
These were also observed in \cite{azulayMultigridAugmentedDeepLearning2023} and other works.
To address (1) one can use multigrid or flexible GMRES \cite{saadFlexibleInnerOuterPreconditioned1993} that in principle allow for nonlinear preconditioners. However, (2) cannot be addressed without excessive training. 

In fact, even a simpler check fails: can the network improve a very coarse solution produced by an approximate CG solver?
Suppose we use Conjugate Gradient (CG) method as the iterative solver. One way to combine the network with CG is to use CG's low-accuracy output as $u_0$ in the network input. The solution $u_0$ is expected to have high-frequency modes more or less already treated.
We find that our network is able to further improve the solution by a limited amount: If CG's output has relative error of 10\%, then our network's output has final error of less than 3\%; If CG's output has relative error of 0.9\%, then our network's output has final error of 0.8\% (train size=500).

\bibliographystyle{elsarticle-num} 
\bibliography{Myref}

@article{azulayMultigridAugmentedDeepLearning2023,
  title = {Multigrid-{{Augmented Deep Learning Preconditioners}} for the {{Helmholtz Equation}}},
  author = {Azulay, Yael and Treister, Eran},
  year = {2023},
  month = jun,
  journal = {SIAM J. Sci. Comput.},
  volume = {45},
  number = {3},
  pages = {S127-S151},
  publisher = {{Society for Industrial and Applied Mathematics}},
  urldate = {2023-09-03}
}

@article{babuskaPollutionEffectFEM1997,
  title = {Is the {{Pollution Effect}} of the {{FEM Avoidable}} for the {{Helmholtz Equation Considering High Wave Numbers}}?},
  author = {Babu{\v s}ka, Ivo M. and Sauter, Stefan A.},
  year = {1997},
  month = dec,
  journal = {SIAM J. Numer. Anal.},
  volume = {34},
  number = {6},
  pages = {2392--2423},
  publisher = {{Society for Industrial and Applied Mathematics}},
  urldate = {2024-01-22}
}

@inproceedings{basriFrequencyBiasNeural2020,
  title = {Frequency {{Bias}} in {{Neural Networks}} for {{Input}} of {{Non-Uniform Density}}},
  booktitle = {Proceedings of the 37th {{International Conference}} on {{Machine Learning}}},
  author = {Basri, Ronen and Galun, Meirav and Geifman, Amnon and Jacobs, David and Kasten, Yoni and Kritchman, Shira},
  year = {2020},
  month = nov,
  pages = {685--694},
  publisher = {PMLR},
  urldate = {2024-04-17},
  langid = {english}
}

@article{bennerSurveyProjectionBasedModel2015,
  title = {A {{Survey}} of {{Projection-Based Model Reduction Methods}} for {{Parametric Dynamical Systems}}},
  author = {Benner, Peter and Gugercin, Serkan and Willcox, Karen},
  year = {2015},
  month = jan,
  journal = {SIAM Rev.},
  volume = {57},
  number = {4},
  pages = {483--531},
  publisher = {{Society for Industrial and Applied Mathematics}},
  urldate = {2024-10-21}
}

@article{bhattacharyaModelReductionNeural2021,
  title = {Model {{Reduction And Neural Networks For Parametric PDEs}}},
  author = {Bhattacharya, Kaushik and Hosseini, Bamdad and Kovachki, Nikola B. and Stuart, Andrew M.},
  year = {2021},
  journal = {The SMAI Journal of computational mathematics},
  volume = {7},
  pages = {121--157},
  urldate = {2023-10-05}
}

@article{caiDeepMMnetInferring2021,
  title = {{{DeepM}}\&{{Mnet}}: {{Inferring}} the Electroconvection Multiphysics Fields Based on Operator Approximation by Neural Networks},
  shorttitle = {{{DeepM}}\&{{Mnet}}},
  author = {Cai, Shengze and Wang, Zhicheng and Lu, Lu and Zaki, Tamer A. and Karniadakis, George Em},
  year = {2021},
  month = jul,
  journal = {Journal of Computational Physics},
  volume = {436},
  pages = {110296},
  urldate = {2024-03-02}
}

@misc{caoChooseTransformerFourier2021,
  title = {Choose a {{Transformer}}: {{Fourier}} or {{Galerkin}}},
  shorttitle = {Choose a {{Transformer}}},
  author = {Cao, Shuhao},
  year = {2021},
  month = nov,
  number = {arXiv:2105.14995},
  eprint = {2105.14995},
  primaryclass = {cs, math},
  publisher = {arXiv},
  urldate = {2024-05-16},
  archiveprefix = {arXiv}
}

@misc{chenMetaMgNetMetaMultigrid2020,
  title = {Meta-{{MgNet}}: {{Meta Multigrid Networks}} for {{Solving Parameterized Partial Differential Equations}}},
  shorttitle = {Meta-{{MgNet}}},
  author = {Chen, Yuyan and Dong, Bin and Xu, Jinchao},
  year = {2020},
  month = nov,
  number = {arXiv:2010.14088},
  eprint = {2010.14088},
  primaryclass = {cs, math},
  publisher = {arXiv},
  urldate = {2024-03-07},
  archiveprefix = {arXiv}
}

@article{chenUniversalApproximationNonlinear1995,
  title = {Universal Approximation to Nonlinear Operators by Neural Networks with Arbitrary Activation Functions and Its Application to Dynamical Systems},
  author = {Chen, T. and Chen, H.},
  year = {1995},
  journal = {IEEE Trans Neural Netw},
  volume = {6},
  number = {4},
  pages = {911--917},
  langid = {english},
  pmid = {18263379}
}

@misc{dosovitskiyImageWorth16x162021,
  title = {An {{Image}} Is {{Worth}} 16x16 {{Words}}: {{Transformers}} for {{Image Recognition}} at {{Scale}}},
  shorttitle = {An {{Image}} Is {{Worth}} 16x16 {{Words}}},
  author = {Dosovitskiy, Alexey and Beyer, Lucas and Kolesnikov, Alexander and Weissenborn, Dirk and Zhai, Xiaohua and Unterthiner, Thomas and Dehghani, Mostafa and Minderer, Matthias and Heigold, Georg and Gelly, Sylvain and Uszkoreit, Jakob and Houlsby, Neil},
  year = {2021},
  month = jun,
  number = {arXiv:2010.11929},
  eprint = {2010.11929},
  primaryclass = {cs},
  publisher = {arXiv},
  urldate = {2023-11-30},
  archiveprefix = {arXiv}
}

@article{eDeepRitzMethod2018,
  title = {The {{Deep Ritz Method}}: {{A Deep Learning-Based Numerical Algorithm}} for {{Solving Variational Problems}}},
  shorttitle = {The {{Deep Ritz Method}}},
  author = {E, Weinan and Yu, Bing},
  year = {2018},
  month = mar,
  journal = {Commun. Math. Stat.},
  volume = {6},
  number = {1},
  pages = {1--12},
  urldate = {2023-08-29},
  langid = {english}
}

@article{eProposalMachineLearning2017,
  title = {A {{Proposal}} on {{Machine Learning}} via {{Dynamical Systems}}},
  author = {E, Weinan},
  year = {2017},
  month = mar,
  journal = {Commun. Math. Stat.},
  volume = {5},
  number = {1},
  pages = {1--11},
  urldate = {2023-08-29},
  langid = {english}
}

@misc{fanaskovSpectralNeuralOperators2024,
  title = {Spectral {{Neural Operators}}},
  author = {Fanaskov, V. and Oseledets, I.},
  year = {2024},
  month = apr,
  number = {arXiv:2205.10573},
  eprint = {2205.10573},
  primaryclass = {cs, math},
  publisher = {arXiv},
  urldate = {2024-09-16},
  archiveprefix = {arXiv}
}

@misc{gopakumarFourierNeuralOperator2023,
  title = {Fourier {{Neural Operator}} for {{Plasma Modelling}}},
  author = {Gopakumar, Vignesh and Pamela, Stanislas and Zanisi, Lorenzo and Li, Zongyi and Anandkumar, Anima and Team, {\relax MAST}},
  year = {2023},
  month = feb,
  number = {arXiv:2302.06542},
  eprint = {2302.06542},
  primaryclass = {physics},
  publisher = {arXiv},
  urldate = {2024-03-13},
  archiveprefix = {arXiv}
}

@inproceedings{guoConvolutionalNeuralNetworks2016,
  title = {Convolutional {{Neural Networks}} for {{Steady Flow Approximation}}},
  booktitle = {Proceedings of the 22nd {{ACM SIGKDD International Conference}} on {{Knowledge Discovery}} and {{Data Mining}}},
  author = {Guo, Xiaoxiao and Li, Wei and Iorio, Francesco},
  year = {2016},
  month = aug,
  series = {{{KDD}} '16},
  pages = {481--490},
  publisher = {Association for Computing Machinery},
  address = {New York, NY, USA},
  urldate = {2023-08-29},
  isbn = {978-1-4503-4232-2}
}

@inproceedings{guptaMultiwaveletbasedOperatorLearning2021,
  title = {Multiwavelet-Based {{Operator Learning}} for {{Differential Equations}}},
  booktitle = {Advances in {{Neural Information Processing Systems}}},
  author = {Gupta, Gaurav and Xiao, Xiongye and Bogdan, Paul},
  year = {2021},
  volume = {34},
  pages = {24048--24062},
  publisher = {Curran Associates, Inc.},
  urldate = {2024-09-18}
}

@inproceedings{heDeepResidualLearning2016,
  title = {Deep {{Residual Learning}} for {{Image Recognition}}},
  booktitle = {2016 {{IEEE Conference}} on {{Computer Vision}} and {{Pattern Recognition}} ({{CVPR}})},
  author = {He, Kaiming and Zhang, Xiangyu and Ren, Shaoqing and Sun, Jian},
  year = {2016},
  month = jun,
  pages = {770--778},
  publisher = {IEEE},
  address = {Las Vegas, NV, USA},
  urldate = {2023-11-21},
  isbn = {978-1-4673-8851-1},
  langid = {english}
}

@article{heMgNetUnifiedFramework2019,
  title = {{{MgNet}}: {{A Unified Framework}} of {{Multigrid}} and {{Convolutional Neural Network}}},
  shorttitle = {{{MgNet}}},
  author = {He, Juncai and Xu, Jinchao},
  year = {2019},
  month = jul,
  journal = {Sci. China Math.},
  volume = {62},
  number = {7},
  eprint = {1901.10415},
  primaryclass = {cs, math},
  pages = {1331--1354},
  urldate = {2024-03-07},
  archiveprefix = {arXiv}
}

@misc{heMgNOEfficientParameterization2023,
  title = {{{MgNO}}: {{Efficient Parameterization}} of {{Linear Operators}} via {{Multigrid}}},
  shorttitle = {{{MgNO}}},
  author = {He, Juncai and Liu, Xinliang and Xu, Jinchao},
  year = {2023},
  month = oct,
  number = {arXiv:2310.19809},
  eprint = {2310.19809},
  primaryclass = {cs, math},
  publisher = {arXiv},
  urldate = {2024-05-14},
  archiveprefix = {arXiv}
}

@misc{ioffeBatchNormalizationAccelerating2015,
  title = {Batch {{Normalization}}: {{Accelerating Deep Network Training}} by {{Reducing Internal Covariate Shift}}},
  shorttitle = {Batch {{Normalization}}},
  author = {Ioffe, Sergey and Szegedy, Christian},
  year = {2015},
  month = mar,
  number = {arXiv:1502.03167},
  eprint = {1502.03167},
  publisher = {arXiv},
  urldate = {2024-10-21},
  archiveprefix = {arXiv}
}

@misc{karrasAliasFreeGenerativeAdversarial2021,
  title = {Alias-{{Free Generative Adversarial Networks}}},
  author = {Karras, Tero and Aittala, Miika and Laine, Samuli and H{\"a}rk{\"o}nen, Erik and Hellsten, Janne and Lehtinen, Jaakko and Aila, Timo},
  year = {2021},
  month = oct,
  number = {arXiv:2106.12423},
  eprint = {2106.12423},
  primaryclass = {cs, stat},
  publisher = {arXiv},
  urldate = {2024-09-20},
  archiveprefix = {arXiv}
}

@article{khooSolvingParametricPDE2021,
  title = {Solving Parametric {{PDE}} Problems with Artificial Neural Networks},
  author = {Khoo, Yuehaw and Lu, Jianfeng and Ying, Lexing},
  year = {2021},
  month = jun,
  journal = {Eur. J. Appl. Math},
  volume = {32},
  number = {3},
  eprint = {1707.03351},
  primaryclass = {math},
  pages = {421--435},
  urldate = {2023-08-29},
  archiveprefix = {arXiv}
}

@article{kissasLearningOperatorsCoupled2022,
  title = {Learning Operators with Coupled Attention},
  author = {Kissas, Georgios and Seidman, Jacob H. and Guilhoto, Leonardo Ferreira and Preciado, Victor M. and Pappas, George J. and Perdikaris, Paris},
  year = {2022},
  month = jan,
  journal = {J. Mach. Learn. Res.},
  volume = {23},
  number = {1},
  pages = {215:9636--215:9698}
}

@article{kontolatiLearningNonlinearOperators2024,
  title = {Learning Nonlinear Operators in Latent Spaces for Real-Time Predictions of Complex Dynamics in Physical Systems},
  author = {Kontolati, Katiana and Goswami, Somdatta and Em Karniadakis, George and Shields, Michael D.},
  year = {2024},
  month = jun,
  journal = {Nat Commun},
  volume = {15},
  number = {1},
  pages = {5101},
  publisher = {Nature Publishing Group},
  urldate = {2024-09-16},
  copyright = {2024 The Author(s)},
  langid = {english}
}

@misc{kopanicakovaDeepOnetBasedPreconditioning2024,
  title = {{{DeepOnet Based Preconditioning Strategies For Solving Parametric Linear Systems}} of {{Equations}}},
  author = {Kopani{\v c}{\'a}kov{\'a}, Alena and Karniadakis, George Em},
  year = {2024},
  month = jan,
  number = {arXiv:2401.02016},
  eprint = {2401.02016},
  primaryclass = {cs, math},
  publisher = {arXiv},
  urldate = {2024-01-16},
  archiveprefix = {arXiv}
}

@article{kovachkiNeuralOperatorLearning,
  title = {Neural {{Operator}}: {{Learning Maps Between Function Spaces With Applications}} to {{PDEs}}},
  author = {Kovachki, Nikola and Li, Zongyi and Liu, Burigede and Azizzadenesheli, Kamyar and Bhattacharya, Kaushik and Stuart, Andrew},
  langid = {english}
}

@misc{lanthalerErrorEstimatesDeepOnets2022,
  title = {Error Estimates for {{DeepOnets}}: {{A}} Deep Learning Framework in Infinite Dimensions},
  shorttitle = {Error Estimates for {{DeepOnets}}},
  author = {Lanthaler, Samuel and Mishra, Siddhartha and Karniadakis, George Em},
  year = {2022},
  month = jan,
  number = {arXiv:2102.09618},
  eprint = {2102.09618},
  primaryclass = {cs, math},
  publisher = {arXiv},
  urldate = {2024-03-02},
  archiveprefix = {arXiv}
}

@misc{lanthalerNonlinearReconstructionOperator2022,
  title = {Nonlinear {{Reconstruction}} for {{Operator Learning}} of {{PDEs}} with {{Discontinuities}}},
  author = {Lanthaler, Samuel and Molinaro, Roberto and Hadorn, Patrik and Mishra, Siddhartha},
  year = {2022},
  month = oct,
  number = {arXiv:2210.01074},
  eprint = {2210.01074},
  primaryclass = {cs, math},
  publisher = {arXiv},
  urldate = {2024-09-16},
  archiveprefix = {arXiv}
}

@article{lecunDeepLearning2015,
  title = {Deep Learning},
  author = {LeCun, Yann and Bengio, Yoshua and Hinton, Geoffrey},
  year = {2015},
  month = may,
  journal = {Nature},
  volume = {521},
  number = {7553},
  pages = {436--444},
  publisher = {Nature Publishing Group},
  urldate = {2024-09-16},
  copyright = {2015 Springer Nature Limited},
  langid = {english}
}

@misc{liFourierNeuralOperator2021,
  title = {Fourier {{Neural Operator}} for {{Parametric Partial Differential Equations}}},
  author = {Li, Zongyi and Kovachki, Nikola and Azizzadenesheli, Kamyar and Liu, Burigede and Bhattacharya, Kaushik and Stuart, Andrew and Anandkumar, Anima},
  year = {2021},
  month = may,
  number = {arXiv:2010.08895},
  eprint = {2010.08895},
  primaryclass = {cs, math},
  publisher = {arXiv},
  urldate = {2023-08-29},
  archiveprefix = {arXiv}
}

@misc{linGreenMultigridNetwork2024,
  title = {Green {{Multigrid Network}}},
  author = {Lin, Ye and Lee, Young Ju and Jia, Jiwei},
  year = {2024},
  month = jul,
  number = {arXiv:2407.03593},
  eprint = {2407.03593},
  primaryclass = {math},
  publisher = {arXiv},
  urldate = {2025-09-24},
  archiveprefix = {arXiv}
}

@article{linOperatorLearningPredicting2021,
  title = {Operator Learning for Predicting Multiscale Bubble Growth Dynamics},
  author = {Lin, Chensen and Li, Zhen and Lu, Lu and Cai, Shengze and Maxey, Martin and Karniadakis, George Em},
  year = {2021},
  month = mar,
  journal = {The Journal of Chemical Physics},
  volume = {154},
  number = {10},
  pages = {104118},
  urldate = {2024-03-13},
  langid = {english}
}

@misc{liTransformerPartialDifferential2023,
  title = {Transformer for {{Partial Differential Equations}}' {{Operator Learning}}},
  author = {Li, Zijie and Meidani, Kazem and Farimani, Amir Barati},
  year = {2023},
  month = apr,
  number = {arXiv:2205.13671},
  eprint = {2205.13671},
  primaryclass = {cs},
  publisher = {arXiv},
  urldate = {2024-09-19},
  archiveprefix = {arXiv}
}

@inproceedings{longPDENetLearningPDEs2018,
  title = {{{PDE-Net}}: {{Learning PDEs}} from {{Data}}},
  shorttitle = {{{PDE-Net}}},
  booktitle = {Proceedings of the 35th {{International Conference}} on {{Machine Learning}}},
  author = {Long, Zichao and Lu, Yiping and Ma, Xianzhong and Dong, Bin},
  year = {2018},
  month = jul,
  pages = {3208--3216},
  publisher = {PMLR},
  urldate = {2023-08-29},
  langid = {english}
}

@article{longPDENetLearningPDEs2019,
  title = {{{PDE-Net}} 2.0: {{Learning PDEs}} from Data with a Numeric-Symbolic Hybrid Deep Network},
  shorttitle = {{{PDE-Net}} 2.0},
  author = {Long, Zichao and Lu, Yiping and Dong, Bin},
  year = {2019},
  month = dec,
  journal = {Journal of Computational Physics},
  volume = {399},
  pages = {108925},
  urldate = {2023-08-29}
}

@article{luComprehensiveFairComparison2022,
  title = {A Comprehensive and Fair Comparison of Two Neural Operators (with Practical Extensions) Based on {{FAIR}} Data},
  author = {Lu, Lu and Meng, Xuhui and Cai, Shengze and Mao, Zhiping and Goswami, Somdatta and Zhang, Zhongqiang and Karniadakis, George Em},
  year = {2022},
  month = apr,
  journal = {Computer Methods in Applied Mechanics and Engineering},
  volume = {393},
  pages = {114778},
  urldate = {2023-12-01}
}

@article{luDeepXDEDeepLearning2021a,
  title = {{{DeepXDE}}: {{A Deep Learning Library}} for {{Solving Differential Equations}}},
  shorttitle = {{{DeepXDE}}},
  author = {Lu, Lu and Meng, Xuhui and Mao, Zhiping and Karniadakis, George Em},
  year = {2021},
  month = jan,
  journal = {SIAM Rev.},
  volume = {63},
  number = {1},
  pages = {208--228},
  publisher = {{Society for Industrial and Applied Mathematics}},
  urldate = {2023-11-18}
}

@article{luLearningNonlinearOperators2021,
  title = {Learning Nonlinear Operators via {{DeepONet}} Based on the Universal Approximation Theorem of Operators},
  author = {Lu, Lu and Jin, Pengzhan and Pang, Guofei and Zhang, Zhongqiang and Karniadakis, George Em},
  year = {2021},
  month = mar,
  journal = {Nat Mach Intell},
  volume = {3},
  number = {3},
  pages = {218--229},
  publisher = {Nature Publishing Group},
  urldate = {2023-08-29},
  copyright = {2021 The Author(s), under exclusive licence to Springer Nature Limited},
  langid = {english}
}

@article{maoPPDONetDeepOperator2023,
  title = {{{PPDONet}}: {{Deep Operator Networks}} for {{Fast Prediction}} of {{Steady-state Solutions}} in {{Disk}}--{{Planet Systems}}},
  shorttitle = {{{PPDONet}}},
  author = {Mao, Shunyuan and Dong, Ruobing and Lu, Lu and Yi, Kwang Moo and Wang, Sifan and Perdikaris, Paris},
  year = {2023},
  month = jun,
  journal = {ApJL},
  volume = {950},
  number = {2},
  pages = {L12},
  publisher = {The American Astronomical Society},
  urldate = {2024-03-13},
  langid = {english}
}

@article{odenaDeconvolutionCheckerboardArtifacts2016,
  title = {Deconvolution and {{Checkerboard Artifacts}}},
  author = {Odena, Augustus and Dumoulin, Vincent and Olah, Chris},
  year = {2016},
  month = oct,
  journal = {Distill},
  volume = {1},
  number = {10},
  pages = {e3},
  urldate = {2024-03-06},
  langid = {english}
}

@article{oleary-roseberryDerivativeInformedNeuralOperator2024,
  title = {Derivative-{{Informed Neural Operator}}: {{An}} Efficient Framework for High-Dimensional Parametric Derivative Learning},
  shorttitle = {Derivative-{{Informed Neural Operator}}},
  author = {{O'Leary-Roseberry}, Thomas and Chen, Peng and Villa, Umberto and Ghattas, Omar},
  year = {2024},
  month = jan,
  journal = {Journal of Computational Physics},
  volume = {496},
  pages = {112555},
  urldate = {2024-05-07}
}

@misc{pathakFourCastNetGlobalDatadriven2022,
  title = {{{FourCastNet}}: {{A Global Data-driven High-resolution Weather Model}} Using {{Adaptive Fourier Neural Operators}}},
  shorttitle = {{{FourCastNet}}},
  author = {Pathak, Jaideep and Subramanian, Shashank and Harrington, Peter and Raja, Sanjeev and Chattopadhyay, Ashesh and Mardani, Morteza and Kurth, Thorsten and Hall, David and Li, Zongyi and Azizzadenesheli, Kamyar and Hassanzadeh, Pedram and Kashinath, Karthik and Anandkumar, Animashree},
  year = {2022},
  month = feb,
  number = {arXiv:2202.11214},
  eprint = {2202.11214},
  primaryclass = {physics},
  publisher = {arXiv},
  urldate = {2024-03-13},
  archiveprefix = {arXiv},
  langid = {english}
}

@inproceedings{rahamanSpectralBiasNeural2019,
  title = {On the {{Spectral Bias}} of {{Neural Networks}}},
  booktitle = {Proceedings of the 36th {{International Conference}} on {{Machine Learning}}},
  author = {Rahaman, Nasim and Baratin, Aristide and Arpit, Devansh and Draxler, Felix and Lin, Min and Hamprecht, Fred and Bengio, Yoshua and Courville, Aaron},
  year = {2019},
  month = may,
  pages = {5301--5310},
  publisher = {PMLR},
  urldate = {2024-04-17},
  langid = {english}
}

@article{rahmanUshapedNeuralOperators2023,
  title = {U-{{NO}}: {{U-shaped Neural Operators}}},
  shorttitle = {U-{{NO}}},
  author = {Rahman, Md Ashiqur and Ross, Zachary E. and Azizzadenesheli, Kamyar},
  year = {2023},
  month = jan,
  journal = {Transactions on Machine Learning Research},
  urldate = {2023-11-21},
  langid = {english}
}

@article{raissiPhysicsinformedNeuralNetworks2019,
  title = {Physics-Informed Neural Networks: {{A}} Deep Learning Framework for Solving Forward and Inverse Problems Involving Nonlinear Partial Differential Equations},
  shorttitle = {Physics-Informed Neural Networks},
  author = {Raissi, M. and Perdikaris, P. and Karniadakis, G. E.},
  year = {2019},
  month = feb,
  journal = {Journal of Computational Physics},
  volume = {378},
  pages = {686--707},
  urldate = {2023-09-08}
}

@article{raonicCONVOLUTIONALNEURALOPERATORS2023,
  title = {{{CONVOLUTIONAL NEURAL OPERATORS}}},
  author = {Raonic, Bogdan and Rohner, Tobias},
  year = {2023},
  langid = {english}
}

@inproceedings{ronnebergerUNetConvolutionalNetworks2015a,
  title = {{U-Net: Convolutional Networks for Biomedical Image Segmentation}},
  shorttitle = {{U-Net}},
  booktitle = {{Medical Image Computing and Computer-Assisted Intervention (MICCAI)}},
  author = {Ronneberger, Olaf and Fischer, Philipp and Brox, Thomas},
  year = {2015},
  urldate = {2023-11-21},
  langid = {ngerman}
}

@article{saadFlexibleInnerOuterPreconditioned1993,
  title = {A {{Flexible Inner-Outer Preconditioned GMRES Algorithm}}},
  author = {Saad, Youcef},
  year = {1993},
  month = mar,
  journal = {SIAM J. Sci. Comput.},
  volume = {14},
  number = {2},
  pages = {461--469},
  publisher = {{Society for Industrial and Applied Mathematics}},
  urldate = {2024-03-04}
}

@misc{sapplDeepLearningPreconditioners2019,
  title = {Deep {{Learning}} of {{Preconditioners}} for {{Conjugate Gradient Solvers}} in {{Urban Water Related Problems}}},
  author = {Sappl, Johannes and Seiler, Laurent and Harders, Matthias and Rauch, Wolfgang},
  year = {2019},
  month = jun,
  number = {arXiv:1906.06925},
  eprint = {1906.06925},
  primaryclass = {cs, math, stat},
  publisher = {arXiv},
  urldate = {2023-09-01},
  archiveprefix = {arXiv}
}

@misc{sitzmannImplicitNeuralRepresentations2020,
  title = {Implicit {{Neural Representations}} with {{Periodic Activation Functions}}},
  author = {Sitzmann, Vincent and Martel, Julien N. P. and Bergman, Alexander W. and Lindell, David B. and Wetzstein, Gordon},
  year = {2020},
  month = jun,
  number = {arXiv:2006.09661},
  eprint = {2006.09661},
  publisher = {arXiv},
  urldate = {2024-11-18},
  archiveprefix = {arXiv}
}

@misc{tancikFourierFeaturesLet2020a,
  title = {Fourier {{Features Let Networks Learn High Frequency Functions}} in {{Low Dimensional Domains}}},
  author = {Tancik, Matthew and Srinivasan, Pratul P. and Mildenhall, Ben and {Fridovich-Keil}, Sara and Raghavan, Nithin and Singhal, Utkarsh and Ramamoorthi, Ravi and Barron, Jonathan T. and Ng, Ren},
  year = {2020},
  month = jun,
  number = {arXiv:2006.10739},
  eprint = {2006.10739},
  primaryclass = {cs},
  publisher = {arXiv},
  urldate = {2023-11-21},
  archiveprefix = {arXiv}
}

@inproceedings{vaswaniAttentionAllYou2017a,
  title = {Attention Is {{All}} You {{Need}}},
  booktitle = {Advances in {{Neural Information Processing Systems}}},
  author = {Vaswani, Ashish and Shazeer, Noam and Parmar, Niki and Uszkoreit, Jakob and Jones, Llion and Gomez, Aidan N and ukasz Kaiser, {\L} and Polosukhin, Illia},
  year = {2017},
  volume = {30},
  publisher = {Curran Associates, Inc.},
  urldate = {2024-09-20}
}

@misc{wenUFNOEnhancedFourier2022,
  title = {U-{{FNO}} -- {{An}} Enhanced {{Fourier}} Neural Operator-Based Deep-Learning Model for Multiphase Flow},
  author = {Wen, Gege and Li, Zongyi and Azizzadenesheli, Kamyar and Anandkumar, Anima and Benson, Sally M.},
  year = {2022},
  month = may,
  number = {arXiv:2109.03697},
  eprint = {2109.03697},
  primaryclass = {physics},
  publisher = {arXiv},
  urldate = {2023-11-17},
  archiveprefix = {arXiv}
}

@article{winovichConvPDEUQConvolutionalNeural2019,
  title = {{{ConvPDE-UQ}}: {{Convolutional}} Neural Networks with Quantified Uncertainty for Heterogeneous Elliptic Partial Differential Equations on Varied Domains},
  shorttitle = {{{ConvPDE-UQ}}},
  author = {Winovich, Nick and Ramani, Karthik and Lin, Guang},
  year = {2019},
  month = oct,
  journal = {Journal of Computational Physics},
  volume = {394},
  pages = {263--279},
  urldate = {2024-11-18}
}

@misc{yePDEformerFoundationModel2024,
  title = {{{PDEformer}}: {{Towards}} a {{Foundation Model}} for {{One-Dimensional Partial Differential Equations}}},
  shorttitle = {{{PDEformer}}},
  author = {Ye, Zhanhong and Huang, Xiang and Chen, Leheng and Liu, Hongsheng and Wang, Zidong and Dong, Bin},
  year = {2024},
  month = apr,
  number = {arXiv:2402.12652},
  eprint = {2402.12652},
  primaryclass = {cs, math},
  publisher = {arXiv},
  urldate = {2024-05-14},
  archiveprefix = {arXiv}
}

@misc{yuMultiScaleContextAggregation2016,
  title = {Multi-{{Scale Context Aggregation}} by {{Dilated Convolutions}}},
  author = {Yu, Fisher and Koltun, Vladlen},
  year = {2016},
  month = apr,
  number = {arXiv:1511.07122},
  eprint = {1511.07122},
  primaryclass = {cs},
  publisher = {arXiv},
  urldate = {2023-11-30},
  archiveprefix = {arXiv}
}

@article{zhangBlendingNeuralOperators2024,
  title = {Blending Neural Operators and Relaxation Methods in {{PDE}} Numerical Solvers},
  author = {Zhang, Enrui and Kahana, Adar and Kopani{\v c}{\'a}kov{\'a}, Alena and Turkel, Eli and Ranade, Rishikesh and Pathak, Jay and Karniadakis, George Em},
  year = {2024},
  month = nov,
  journal = {Nat Mach Intell},
  volume = {6},
  number = {11},
  pages = {1303--1313},
  publisher = {Nature Publishing Group},
  urldate = {2025-09-24},
  copyright = {2024 The Author(s), under exclusive licence to Springer Nature Limited},
  langid = {english}
}

@article{zhuBayesianDeepConvolutional2018,
  title = {Bayesian Deep Convolutional Encoder--Decoder Networks for Surrogate Modeling and Uncertainty Quantification},
  author = {Zhu, Yinhao and Zabaras, nicholas},
  year = {2018},
  month = aug,
  journal = {Journal of Computational Physics},
  volume = {366},
  pages = {415--447},
  publisher = {Academic Press},
  urldate = {2023-08-29},
  langid = {american}
}

@misc{zhangHybridIterativeNumerical2022,
  title = {A {{Hybrid Iterative Numerical Transferable Solver}} ({{HINTS}}) for {{PDEs Based}} on {{Deep Operator Network}} and {{Relaxation Methods}}},
  author = {Zhang, Enrui and Kahana, Adar and Turkel, Eli and Ranade, Rishikesh and Pathak, Jay and Karniadakis, George Em},
  year = {2022},
  month = aug,
  number = {arXiv:2208.13273},
  eprint = {2208.13273},
  primaryclass = {cs, math},
  publisher = {arXiv},
  urldate = {2024-01-16},
  archiveprefix = {arXiv}
}

@misc{hu_hybrid_2024,
    title = {A hybrid iterative method based on {MIONet} for {PDEs}: {Theory} and numerical examples},
    shorttitle = {A hybrid iterative method based on {MIONet} for {PDEs}},
    url = {http://arxiv.org/abs/2402.07156},
    doi = {10.48550/arXiv.2402.07156},
    urldate = {2025-09-24},
    publisher = {arXiv},
    author = {Hu, Jun and Jin, Pengzhan},
    month = feb,
    year = {2024},
    note = {arXiv:2402.07156 [math]},
}

@book{quarteroni2015reduced,
  title={Reduced basis methods for partial differential equations: an introduction},
  author={Quarteroni, Alfio and Manzoni, Andrea and Negri, Federico},
  volume={92},
  year={2015},
  publisher={Springer}
}

@article{khoo2019switchnet,
  title={{SwitchNet}: a neural network model for forward and inverse scattering problems},
  author={Khoo, Yuehaw and Ying, Lexing},
  journal={SIAM Journal on Scientific Computing},
  volume={41},
  number={5},
  pages={A3182--A3201},
  year={2019},
  publisher={SIAM}
}


%



\end{document}